\documentclass[3p]{elsarticle}

\usepackage{endnotes}
\usepackage{soul}
\usepackage{booktabs}
\usepackage{longtable}
\usepackage{amsfonts}
\usepackage{amsmath}
\usepackage{amssymb}

\usepackage{geometry}
\usepackage{setspace}
\usepackage{caption}
\usepackage{enumitem}
\usepackage{algorithm}
\usepackage{algpseudocode}

\usepackage{natbib}
\bibpunct[, ]{(}{)}{,}{a}{}{,}%
\def\bibfont{\small}%

\usepackage{makecell}

\newtheorem{definition}{Definition}
\newtheorem{remark}{Remark}

\biboptions{authoryear}

\makeatletter
\def\ps@pprintTitle{%
  \let\@oddhead\@empty
  \let\@evenhead\@empty
  \let\@oddfoot\@empty
  \let\@evenfoot\@oddfoot
}
\makeatother

\begin{document}

\begin{frontmatter}

\title{Vehicle Routing with Stochastic Demands and\\ Partial Reoptimization}

\author[poly]{Alexandre M. Florio\corref{mycorrespondingauthor}}
\cortext[mycorrespondingauthor]{Corresponding author}
\ead{aflorio@gmail.com}

\author[emse]{Dominique Feillet}
\ead{feillet@emse.fr}

\author[puc]{Marcus Poggi}
\ead{poggi@inf.puc-rio.br}

\author[poly,puc]{Thibaut Vidal}
\ead{thibaut.vidal@polymtl.ca}

\address[poly]{CIRRELT \& SCALE-AI Chair in Data-Driven Supply Chains, Department of Mathematics and Industrial Engineering, Polytechnique Montr\'eal, Canada}
\address[emse]{Mines Saint-\'{E}tienne, Univ. Clermont Auvergne, INP Clermont Auvergne, CNRS, UMR 6158 LIMOS, Saint-\'{E}tienne, France}
%\address[emse]{Mines Saint-\'{E}tienne, Univ. Clermont Auvergne, CNRS, UMR 6158 LIMOS, Centre CMP, Gardanne, France}
%\address[emse]{\'{E}cole des Mines de Saint-\'{E}tienne, Gardanne, France}
\address[puc]{Pontif\'{i}cia Universidade Cat\'{o}lica, Rio de Janeiro, Brazil}

\begin{abstract}
We consider the vehicle routing problem with stochastic demands (VRPSD), a problem in which customer demands are known in distribution at the route planning stage and revealed during route execution upon arrival at each customer. A long-standing open question on the VRPSD concerns the benefits of allowing, during route execution, partial reordering of the planned customer visits. Given the practical importance of this question and the growing interest on the VRPSD under optimal restocking, we study the VRPSD under a recourse policy known as the switch policy. The switch policy is a canonical reoptimization policy that permits only pairs of successive customers to be reordered. We consider this policy jointly with optimal preventive restocking and introduce a branch-cut-and-price algorithm to compute optimal a priori routing plans in this context. At its core, this algorithm features pricing routines where value functions represent the expected cost-to-go along planned routes for all possible states and reordering decisions. To ensure pricing tractability, we adopt a strategy that combines elementary pricing with completion bounds of varying complexity, and solve the pricing problem without relying on dominance rules. Our numerical experiments demonstrate the effectiveness of the algorithm for solving instances with up to 50 customers. Notably, they also give us new insights into the value of reoptimization. The switch policy enables significant cost savings over optimal restocking when the planned routes come from an algorithm built on a deterministic approximation of the data, an important scenario given the difficulty of finding optimal VRPSD solutions. The benefits are smaller when comparing optimal a priori VRPSD solutions obtained for both recourse policies. As it appears, further cost savings may require joint reordering and reassignment of customer visits among vehicles when the context permits.
\end{abstract}


\end{frontmatter}

\section{Introduction}
The vehicle routing problem with stochastic demands (VRPSD) is currently one of the most studied stochastic variants of the VRP. In the VRPSD, customer demands are only known in probability distribution at the route planning stage. Demands are revealed sequentially during route execution, upon arrival of the vehicle at each customer. The problem consists in planning a set of routes such that each customer is visited by one vehicle and the total expected routing cost is minimized. The VRPSD is a canonical problem in transport logistics with a variety of applications including courier services, waste collection, sludge disposal, emergency logistics, delivery of heating oil and industrial gases, replenishment of vending machines and cash logistics \citep[see e.g.][]{lambert1993designing,chan2001multiple,singer2002fleet,ChepuriHomemDeMello2005,markov2020waste}.

We consider the VRPSD under the a priori optimization paradigm, in which the problem is solved for a given recourse policy. In the VRPSD, the actual realized demand along a route may exceed the vehicle capacity. When customer demands must be fully served, recourse policies define a set of rules and actions to guarantee capacity feasibility during route execution. Typically, these rules prescribe replenishment trips to the depot, and in these cases the recourse corresponds to a restocking policy. Traditional restocking policies include the detour-to-depot, or classical, policy \citep{Droretal1989}, and the optimal restocking policy \citep{YeeGolden1980}. While this work focuses on the a priori VRPSD under restocking recourse, there are other modeling approaches for the problem, such as non-restocking recourse policies \citep[e.g., the vehicle pairing strategy by][]{AkErera2007}, and methods based on chance-constrained programming \citep{Droretal1993,noorizadegan2018vehicle,Dinh_2018} and robust optimization \citep{Gounaris_2013,ghosal2020distributionally}. We refer to \cite{Gendreauetal2016} and \cite{oyola2017stochastic,oyola2018stochastic} for more extensive reviews on the VRPSD.

The first generation of algorithms for the VRPSD solved the problem under detour-to-depot restocking. This policy is \emph{reactive}, since it allows restocking only when there is not enough remaining load to serve the current customer, an event known as \emph{failure}. Under classical recourse, the state-of-the-art VRPSD algorithms are the branch-cut-and-price from \cite{Gauvinetal2014}, and the integer L-shaped method from \cite{Jabalietal2014}. More recently, a new generation of algorithms for the VRPSD under optimal restocking have been proposed. These include the integer L-shaped methods from \cite{Louveaux_2018} and \cite{salavati2019exact}, and the branch-cut-and-price from \cite{Florio_2020}. The optimal restocking policy is \emph{preventive}, and may prescribe restocking trips in anticipation of failures. In theory, cost savings when implementing optimal restocking may be significant \citep{Bertazzi_2020}, but, in practice, savings are only marginal when routes serve an expected demand not larger than the capacity of the vehicle \citep{florio2020optimal}.

As a general trend, research on the VRPSD has progressed towards more sophisticated recourse policies as a way to achieve cost savings \citep{Gendreauetal2016}. Yet, whereas recourse policies based on optimized restocking trips have been comprehensively evaluated, policies that involve possible customer visiting order reoptimization have been notably less studied. This is likely due to the difficulty of solving the VRPSD to optimality in such cases. Thus, even if recourse policies based on visiting order reoptimization appear as intuitive to promote a better utilization of the residual vehicle capacity before replenishment trips (e.g., when customers located in close proximity have large differences of demand volatility), it remains a long-standing open research question to characterize their potential benefits.

To progress on this challenging research question, this paper conducts a comprehensive study of a restocking policy called \emph{switch policy}, which allows the actual visiting sequence to differ from the planned sequence by swapping the order of any two customers planned to be visited consecutively. The switch policy is perhaps the simplest policy where changing the visiting order defined at the planning stage is permitted. Arguably, if reoptimization leads to sizable gains, then one can expect to realize at least a share of those gains with such a simple policy.

Models for the VRPSD in which the customer visiting sequence is entirely decided during route execution are proposed under the so-called reoptimization paradigm, where no a priori plan is defined. However, even single-vehicle reoptimization models are intractable in instances with more than 15 customers \citep{SecomandiMargot2009}, and to date can only be solved heuristically with methods from approximate dynamic programming \citep{secomandi2000comparing,Secomandi2001,NovoaStorer2009}. Partial reoptimization policies like the switch policy can be computed for a fixed planned sequence, and they can be used to approximate the cost of full reoptimization models. Finding the optimal sequence under such a policy is still an open problem though \citep{SecomandiMargot2009}. This paper proposes an algorithm to find optimal multi-vehicle routing plans to be executed under the switch policy. Our work bridges several research gaps identified by \cite{SecomandiMargot2009}, namely, the optimization of planned sequences under a reoptimization recourse policy, the measure of the benefit of reoptimization, and the computation of optimal multi-vehicle plans under reoptimization.

To find optimal planned sequences under the switch policy, we propose a novel branch-cut-and-price algorithm. This algorithm relies on a backward labeling procedure for solving the pricing problem, in which labels store value functions that represent the expected cost-to-go along planned routes for all possible states and reordering decisions. A major obstacle when solving the VRPSD with branch-and-price is the derivation of effective dominance rules to discard unpromising partial paths in the pricing algorithm \citep{Florio_2020}. To circumvent this difficulty, our algorithm features an elementary pricing strategy where strong completion bounds are employed to discard unpromising paths. The algorithm is further enhanced with the separation of rounded capacity and subset row inequalities.

In summary, the main contributions of this paper are the following:
\begin{enumerate}[label={(\arabic*)}]
\item In order to assess the benefit of reoptimizing a priori sequences, we introduce the switch policy as a new recourse policy for the VRPSD. We propose a dynamic programming procedure for computing optimal decisions and the expected cost of a customer sequence executed under the switch policy, and a branch-cut-and-price algorithm for finding optimal planned sequences under this recourse policy.
\item The design of such an algorithm offers insights on how to adapt the branch-and-price framework to handle two-stage stochastic problems with complex recourse. When the recourse policy is characterized by rich action and state spaces, dominance between partial routes must occur over a high-dimensional resource space and may not be effective. To improve tractability, we propose a new solution strategy that combines elementary route pricing with strong completion bounds.
\item We measure and report, for the first time, the benefit of partial reoptimization in optimal and non-optimal a priori sequences. More specifically, we compute the value of the switch recourse policy relative to optimal restocking in optimal and non-optimal routing plans, under varying load factors and demand variability settings.
\end{enumerate}

The remainder of this paper is organized as follows: in Section \ref{sec:probdef}, we propose a classification scheme for restocking-based recourse policies, formally define the VRPSD and present mathematical formulations for the problem. In Section \ref{sec:switchpol}, we present a stochastic dynamic programming algorithm for computing the expected cost and optimal restocking and reordering decisions when executing a route under the switch policy. In Section \ref{sec:bcnpalg}, we introduce an exact algorithm for the VRPSD under the switch policy. We discuss the results of extensive numerical experiments in Section \ref{sec:results}. Finally, in Section \ref{sec:conclusions} we provide concluding remarks and perspectives for future research.

\section{Problem Definition} \label{sec:probdef}
The VRPSD is defined on a complete graph $G=(V,A)$, where $V=\{0,\ldots,N\}$ is the set of nodes and $A=\{(i,j):i,j\in V\}$ is the set of arcs. Node $0$ identifies the depot and nodes $1,\ldots,N$ represent customers. A cost $c_{ij}$ is associated to each arc $(i,j)\in A$, and we assume that the triangle inequality holds. A fleet of $K$ vehicles is available at the depot, each with a capacity of $Q$. Customer demands are stochastic with known probability distributions. We assume that demands are discrete and statistically independent. The demand of customer $i$ is represented by the random variable $D_{i}$, and the probability mass function of $D_{i}$ is denoted by $p_{D_{i}}(\cdot)$. Without loss of generality, we assume a \emph{delivery} operation, that is, demands correspond to quantities of a certain commodity that must be delivered to each customer.

Planned routes are defined as follows:
\vspace{3pt}
\begin{definition}[Planned route]
A planned route is a sequence of customers $(s_{1},\ldots,s_{H})$ such that $s_{i}\neq s_{j}$ if $i\neq j$.
\end{definition}
\vspace{3pt}

In this paper, we focus on restocking-based recourse policies (as opposed to, e.g., penalty-based recourse), which guarantee that demands are fully served. A recourse policy defines the set of rules that the driver follows when executing a planned route. These rules may, e.g., instruct the driver to replenish the vehicle or to change the customer visiting order. Concerning restocking decisions, recourse policies may be reactive (replenishment allowed only upon failure) or preventive (replenishment also allowed between customer visits). Next, we propose a classification of preventive recourse policies.
\vspace{3pt}
\begin{definition}[$\pi(k)$ and $\pi^{*}(k)$ policies]
With respect to a planned route $(s_{1},\ldots,s_{H})$, a recourse policy $\pi(k)$, $k\in\mathbb{N}$, is a policy that allows preventive restocking and also allows a customer $s_{h}$ to be visited in a new position $h^{\prime}\in\{\max\{1,h-k\},\ldots,\min\{H,h+k\}\}$. We refer to the optimal policy within a family as $\pi^{*}$.
\end{definition}
\vspace{3pt}

With these definitions, $\pi(0)$ corresponds to a non-optimal policy that allows preventive replenishment, such as the rule-based policy employed by \cite{salavati2019rule}; $\pi^{*}(0)$ refers to the optimal restocking policy as used in the current state-of-the-art algorithms for the VRPSD, and $\pi^{*}(1)$ is the switch policy studied in this paper. Further, $\pi^{*}(H)$ corresponds to fully reoptimizing the route, that is, a policy in which the customer visiting sequence is completely defined during route execution. In this case, the concept of planned route does not exist anymore, as a priori decisions only relate to the assignment of customers to vehicles.

A posteriori routes are defined as follows:
\vspace{3pt}
\begin{definition}[A posteriori route]
For a planned route $(s_{1},\ldots,s_{H})$ and a recourse policy $\pi(k)$, an a posteriori route is a sequence $(u_{1},\ldots,u_{H})$ such that $u_{i}\neq u_{j}$ if $i\neq j$, and $u_{h}\in\{s_{l}:l\geq\max(h-k,1)\wedge l\leq\min(h+k,H)\}$ for all $h\in\{1,\ldots,H\}$.
\end{definition}
\vspace{3pt}

The complexity increases significantly when moving from a policy $\pi(0)$ to a policy $\pi(1)$. There is only one possible a posteriori route when executing a route under policy $\pi(0)$. On the other hand, a planned route with as few as six customers executed under policy $\pi(1)$ may lead to 13 distinct a posteriori routes. In fact, as one easily verifies by taking customers along a planned route $(s_{1},\ldots,s_{H})$ pairwise, the number of a posteriori routes under policy $\pi(1)$ grows exponentially with~$H/2$.

\subsection{Arc-based Formulation} \label{secarcbased}
Assuming a decision vector $\mathbf{x}=\{x_{ij}:i,j\in\{0,\ldots,N\}\}$ indicating the arcs used in the solution, a two-index formulation of the VRPSD under recourse policy $\pi$ is given by:
\begin{align}
&\text{minimize}&&\sum_{i=0}^{N}\sum_{i=0}^{N}c_{ij}x_{ij}+\mathbb{E}\left[\mathcal{Q}^{\pi}(\mathbf{x})\right], \nonumber \\*
&\text{subject to}&&\sum_{j=1}^{N}x_{0j}\leq K,& \nonumber \\*
&&&\sum_{i=0}^{N}x_{ij}=1,&j\in\{1,\ldots,N\}, \nonumber \\*
&&&\sum_{j=0}^{N}x_{ij}=1,&i\in\{1,\ldots,N\}, \nonumber \\*
&&&\sum_{i\in S}\sum_{j\notin S}x_{ij}\geq\left\lceil\sum_{i\in S}\mathbb{E}[D_{i}]/fQ\right\rceil,&S\subset\{1,\ldots,N\}, \label{rcc} \\*
&&&x_{ij}\in\{0,1\}, \nonumber
\end{align}
where $f$ is the load factor parameter that enables (when $f>1$) routes to serve a total expected demand above the vehicle capacity. Most exact methods for the VRPSD assume $f=1$. However, allowing larger load factors may lead to superior solutions with regard to both cost and number of required vehicles, and differences between restocking policies become significant only when $f>1$ \citep{florio2021branch,florio2020optimal}.

The random variable $\mathcal{Q}^{\pi}(\mathbf{x})$ denotes the recourse cost associated with first stage decision $\mathbf{x}$ under policy $\pi$. Even though the arc-based formulation is sometimes interpreted as a two-stage stochastic programming model, the restocking-based VRPSD is actually a multistage problem since, in addition to an initial decision concerning planned routes, decisions are also made after every customer is served \citep{Droretal1989}. The random variable $\mathcal{Q}^{\pi}(\mathbf{x})$ corresponds to the cost of all stages except the first one. State-of-the-art integer L-shaped algorithms for the VRPSD use this formulation as a starting point. In our branch-cut-and-price algorithm, branching is performed on variables of the two-index formulation, and inequalities \eqref{rcc} are separated and lifted to improve the linear bound of the path-based formulation discussed in the next section.

When only preventive restocking is allowed, the recourse cost $\mathcal{Q}^{\pi}(\mathbf{x})$ is nonnegative. When customer reordering is allowed, however, $\mathcal{Q}^{\pi}(\mathbf{x})$ may be negative. This may pose additional challenges when attempting to solve the arc-based formulation with the integer L-shaped method, since the first stage decision does not represent a lower bound on the total solution cost anymore.

\subsection{Path-based Formulation}
The cost of a planned route $r$ under a policy $\pi$ is denoted by a random variable $C_{r}^{\pi}$. Note that $C_{r}^{\pi}$ includes both the planned and the recourse costs of a route that starts at the depot, follows the customer sequence $r$, and finishes at the depot. Given decision variables $\lambda_{r}$ to identify whether route $r$ belongs to the solution, the VRPSD can be modeled as a set-partitioning problem:
\begin{align}
&\text{minimize}&&\sum_{r\in\mathcal{R}}\mathbb{E}[C_{r}^{\pi}]\lambda_{r}, \nonumber \\
&\text{subject to}&&\sum_{r\in\mathcal{R}}\lambda_{r}\leq K, \label{numVctr} \\
&&&\sum_{r\in\mathcal{R}}[i\in r]\lambda_{r}=1,&i\in\{1,\ldots,N\}, \label{spc} \\
&&&\lambda_{r}\in\{0,1\},&r\in\mathcal{R}, \nonumber
\end{align}
where $\mathcal{R}$ is the set of all feasible planned routes. The square brackets are Iverson brackets (i.e., $[P]=1$ if the logical proposition $P$ is true, and $[P]=0$ otherwise). The formulation is further strengthened with subset row cuts defined over customer triplets \citep{jepsen2008subset}:
\begin{align}
\sum_{r\in\mathcal{R}}\left\lfloor\frac{1}{2}\sum_{i\in S}[i\in r]\right\rfloor\lambda_{r}&\leq1,&S\subset\{1,\ldots,N\},|S|=3,\label{eqsrc}\\
\intertext{and with strengthened capacity cuts \citep{baldacci2008exact}:}
\sum_{r\in\mathcal{R}}[r\cap S\neq\emptyset]\lambda_{r}&\geq\left\lceil\sum_{i\in S}\mathbb{E}[D_{i}]/fQ\right\rceil,&S\subset\{1,\ldots,N\}.\label{eqscc}
\end{align}

\section{The $\pi(1)$ (Switch) Policy} \label{sec:switchpol}
The switch policy is a special case of the sliding heuristic (SH) by \cite{SecomandiMargot2009}. Whereas in that work the policy is solved by restricting the state space of a (full) reoptimization model, in this section we present a specialized dynamic programming algorithm for computing the expected cost and the optimal decisions of a planned route executed under the switch policy.

Consider a planned route $r=(s_{1},\ldots,s_{H})$. When the vehicle performs $r$ under the switch policy, we represent the execution state just after a customer is fully served by a triplet $(h,n,q)$, where $h\in\{1,\ldots,H\}$ is the number of customers already served, $n\in\{s_{l}:l\geq\max(h-1,1)\wedge l\leq\min(h+1,H)\}$ is the last customer served, and $q$ is the remaining load in the vehicle. We denote by $\overline{\nu}(h,n,q)$ the expected remaining cost along $r$ when at state $(h,n,q)$.

To facilitate the presentation of the dynamic programming algorithm for the switch policy, we define the functionals $\phi'(i,j,q,\nu(\cdot))$ and $\phi''(i,j,\nu(\cdot))$, which correspond to the expected costs related to the decisions of visiting customer $j$ from node $i$ directly (given a remaining vehicle load of $q$) or after capacity replenishment, respectively, assuming a cost-to-go function $\nu:\{0,\ldots,Q\}\mapsto\mathbb{R}$ once customer $j$ is served:
\begin{equation}\label{eqdefpip}
\phi'(i,j,q,\nu(\cdot))=c_{ij}+\sum_{d=0}^{\infty}[\Gamma_{d,q}(c_{j0}+c_{0j})+\nu(q+Q\Gamma_{d,q}-d)]p_{D_{j}}(d),
\end{equation}
\begin{equation}\label{eqdefpipp}
\phi''(i,j,\nu(\cdot))=c_{i0}+c_{0j}+\sum_{d=0}^{\infty}[\Gamma_{d,Q}(c_{j0}+c_{0j})+\nu(Q+Q\Gamma_{d,Q}-d)]p_{D_{j}}(d),
\end{equation}
where $\Gamma_{d,q}=\max\{0, \lceil(d-q)/Q\rceil\}$ is the number of return trips to the depot that are needed to fully serve a realized demand of $d$ when the remaining load is $q$. We also define $\phi^{*}(i,j,q,\nu(\cdot))=\min\{\phi'(i,j,q,\nu(\cdot)),\phi''(i,j,\nu(\cdot))\}$ as the expected remaining cost given that the best decision between a direct visit or a visit after replenishment is chosen.

\vspace{3pt}
\begin{remark} \label{remarkpenalty}
In some VRPSD models \citep[e.g.,][]{Louveaux_2018}, a failure while serving a customer is penalized with a fixed cost. The idea is to discourage split deliveries, since they involve customer inconvenience and extra costs such as additional unloading time. The dynamic programming model can be immediately extended to this case, by adding extra terms in \eqref{eqdefpip} and \eqref{eqdefpipp}. For the sake of simplicity, we proceed without considering those penalties. Still, in Section \ref{sec:exppolcomp} we also report results assuming that failures are penalized.
\end{remark}
\vspace{3pt}

The base case of the dynamic programming algorithm is given by:
\begin{equation}\label{eqswitchbasecase}
\overline{\nu}(H,n,q)=c_{n0}.
\end{equation}

The available decisions at other states $(h,n,q)$, $1\leq h<H$, depend on the last visited customer. In general, there are three possibilities. If $n=s_{h}$, then customer $s_{h}$ is visited in its original position in the planned route. So, it is now possible to visit customer $s_{h+1}$ or $s_{h+2}$ (only $s_{h+1}$ if $h=H-1$), directly or after restocking. If $n=s_{h-1}$ (only possible if $h>1$), then customers $s_{h}$ and $s_{h-1}$ are switched in the a posteriori route. Both have already been visited though, so that now it is also possible to visit customer $s_{h+1}$ or $s_{h+2}$ (again, only $s_{h+1}$ if $h=H-1$). The last possibility is $n=s_{h+1}$, in which case customers $s_{h}$ and $s_{h+1}$ are switched, but $s_{h}$ has not been visited yet. So, $s_{h}$ must be visited after $s_{h+1}$. Therefore, the remaining cost in each case is given by:
\begin{equation} \label{eqswitchgen}
\overline{\nu}(h,n,q)=
\begin{cases}
\phi^{*}(n,s_{h},q,\nu_{s_{h}}^{h+1}(\cdot)),& \text{if }n=s_{h+1},\\[1ex]
\phi^{*}(n,s_{H},q,\nu_{s_{H}}^{H}(\cdot)),& \text{if }h=H-1,\\[1ex]
\min\{\phi^{*}(n,s_{h+1},q,\nu_{s_{h+1}}^{h+1}(\cdot)),\\[0.5ex] \phantom{\min\{}\phi^{*}(n,s_{h+2},q,\nu_{s_{h+2}}^{h+1}(\cdot))\},&\text{otherwise},
\end{cases}
\end{equation}
where $\nu_{n}^{h}:\{0,\ldots,Q\}\mapsto\mathbb{R}$ is the function defined by $\nu_{n}^{h}(q)=\overline{\nu}(h,n,q)$.

When the route starts, the vehicle proceeds fully loaded from the depot to either the first or the second customer. So, the expected cost of $r$ under the switch policy is given by:
\begin{equation} \label{eq:expcost}
\mathbb{E}[C_{r}^{\pi^{*}\!(1)}]=\min\{\phi'(0,s_{1},Q,\nu_{s_{1}}^{1}(\cdot)),\phi'(0,s_{2},Q,\nu_{s_{2}}^{1}(\cdot))\},
\end{equation}
and can be calculated by solving recursion \eqref{eqswitchgen} by backward induction, starting with the base case~\eqref{eqswitchbasecase}. Note that the dynamic programming algorithm outputs not only \eqref{eq:expcost}, but also a function mapping each state to its associated optimal decision. Concerning complexity, assuming that the cardinalities of the support sets of the demand distributions $p_{D_{i}}(\cdot)$, $i\in\{1,\ldots,N\}$, are $\mathcal{O}(Q)$, the dynamic programming algorithm runs in $\mathcal{O}(Q^2N)$ time.

\section{Branch-Cut-and-Price Algorithm} \label{sec:bcnpalg}
Branch-cut-and-price is a leading method for solving many VRP variants that can be formulated as set-partitioning problems. The method consists in applying branch-and-cut to the set-partitioning formulation of the problem, where, due to the very large number of variables in that formulation, linear bounds are computed by column generation. In most cases, the subproblem (or pricing problem), which is the problem of identifying profitable variables, is a variant of the elementary resource-constrained shortest-path problem \citep{feillet2004exact}. By and large, modern algorithms relax the route elementarity condition and solve the subproblem by a labeling procedure that relies on dominance rules for tractability. We refer to \cite{costa2019exact} for a review of the framework.

\subsection{Overview}
Two main challenges must be overcome when applying the branch-cut-and-price framework to the VRPSD under policy $\pi^{*}(1)$. First, the expected cost of a partial planned route is not precisely defined, since it depends on the previous customer (assuming a partial route from a customer up to the depot), which is unknown. We manage this issue by storing value functions in the labels of the pricing algorithm. These functions represent the expected cost-to-go along the planned route for all possible states and reordering decisions. The second challenge relates to the tractability of the pricing algorithm. When labels store complete value functions, label dominance is hindered because it requires value function dominance over all states. To ensure pricing tractability, we adopt a new strategy that combines elementary pricing with strong completion bounds and solve the pricing problem without relying on dominance rules. In the remainder of this section, we detail the key components of our branch-cut-and-price algorithm, emphasizing the new methodological contributions: the labeling procedure based on value functions (Section \ref{sec:labproc}) and the completion bounds (Section \ref{sec:cbounds}).

\subsection{Pricing Problem}
We refer to the continuous relaxation of the set-partitioning formulation restricted to a subset of the variables as the restricted master problem (RMP). Within column generation, the pricing problem consists in identifying profitable variables to add to the RMP. We denote by $\mathcal{J}$ and $\mathcal{C}$ the sets containing customer triplets and customer sets, respectively, for which subset row cuts \eqref{eqsrc} and strengthened capacity cuts \eqref{eqscc} have been separated and added to the RMP. Given a solution to the RMP, the reduced cost of a route $r$ under policy $\pi^{*}(1)$, denoted by $\overline{C}_{r}$ (for notational simplicity, from now on we omit the superscript identifying the policy), is given by:
\begin{equation} \label{eqredcost}
\overline{C}_{r}=\mathbb{E}[C_{r}]-\kappa-\sum_{i\in r}\alpha_{i}-\sum_{S\in\mathcal{J}}\left\lfloor\sum_{i\in S}[i\in r]/2\right\rfloor\beta_{S}-\sum_{S\in\mathcal{C}}[S\cap r\neq\emptyset]\gamma_{S},
\end{equation}
where $\kappa$, $\alpha_{i}$ ($i\in\{1,\ldots,N\}$), $\beta_{S}$ ($S\in\mathcal{J}$), and $\gamma_{S}$ ($S\in\mathcal{C}$) are the dual values in the solution to the RMP associated with constraints \eqref{numVctr}, \eqref{spc}, \eqref{eqsrc} and \eqref{eqscc}, respectively. Hence, the pricing problem calls for routes $r\in\mathcal{R}$ such that $\overline{C}_{r}<0$.

\subsection{Labeling Procedure} \label{sec:labproc}
The pricing problem is solved with a backward labeling algorithm, where each label stores information about a partial planned route from a customer up to the depot. More specifically, a label $\mathcal{L}$ stores a reference $\mathcal{L}.\mathcal{P}$ to the predecessor label, an ordered-set $\mathcal{L}.\delta$ representing the partial \emph{planned} route, of which the current and next customers are referred to as $\mathcal{L}.n$ and $\mathcal{L}.n'$, respectively. To keep track of the cost-to-go without knowledge about the previous planned node, a label also stores the value functions $\mathcal{L}.\nu,\mathcal{L}.\nu':\{0,\ldots,Q\}\mapsto\mathbb{R}$, which map, for all possible remaining load values, the expected remaining cost along $\mathcal{L}.\delta$ when customer $\mathcal{L}.n$ is served in its planned position and when customers $\mathcal{L}.n$ and $\mathcal{L}.n'$ are switched in the a posteriori route, respectively. Table \ref{tab:label} summarizes all attributes of a label along with their initialization values.

\begin{table}[t]
\centering
\begingroup
\small
\renewcommand{\arraystretch}{0.825}
\caption{\label{tab:label}Label Attributes and Initialization Values}
\begin{tabular}{llr}
\toprule
Attribute & Description & Initialization$^a$ \\
\midrule
$\mathcal{P}$ & Reference to predecessor label. & $\varnothing$ \\
$\delta$ & Partial planned route (sequence). & $(j)$ \\
$n$ \emph{(shortcut)} & Current customer in the planned route. & $j$ \\
$n'$ \emph{(shortcut)} & Next customer in the planned route. & $\varnothing$ \\
$\nu(q)$ & \makecell[tl]{Expected remaining cost just after customer $n$\\is served in its planned position.} & $c_{j0},\quad\forall q\in\{0,\ldots,Q\}$ \\
$\nu'(q)$ & \makecell[tl]{Expected remaining cost just after customer $n'$\\is served, considering that $n$ and $n'$ are switched\\in the a posteriori route.} & $\infty,\quad\forall q\in\{0,\ldots,Q\}$ \\
$\Theta$ & Sum of dual values associated with $\delta$. & $\kappa+\alpha_{j}+\sum_{S\in\mathcal{C}}[j\in S]\gamma_{S}$\\
\bottomrule
\end{tabular}
\captionsetup{justification=centering}
\caption*{$^a$Assuming initialization to a customer $j$.}
\endgroup
\end{table}

Consider a label $\mathcal{L}$. Before $\overline{C}_{\mathcal{L}.\delta}$ can be computed, $\mathbb{E}[C_{\mathcal{L}.\delta}]$ must be determined. Taking into account the information already stored in $\mathcal{L}$, $\mathbb{E}[C_{\mathcal{L}.\delta}]$ is given by:
\begin{equation} \nonumber
\mathbb{E}[C_{\mathcal{L}.\delta}]=
\begin{cases}
\phi'(0,\mathcal{L}.n,Q,\mathcal{L}.\nu),&\text{if }\mathcal{L}.\mathcal{P}=\varnothing,\\
\min\{\phi'(0,\mathcal{L}.n,Q,\mathcal{L}.\nu), \phi'(0,\mathcal{L}.n',Q,\mathcal{L}.\nu')\},&\text{otherwise}.
\end{cases}
\end{equation}

An extension label $\mathcal{E}$ is created by extending $\mathcal{L}$ to a customer $i\notin\mathcal{L}.\delta$. The attributes of $\mathcal{E}$ are defined as follows:
\begin{align}
\mathcal{E}.\mathcal{P}&=\mathcal{L},\nonumber \\
\mathcal{E}.\delta&=(i)\oplus\mathcal{L}.\delta,\nonumber \\
\mathcal{E}.n&=i,\nonumber \\
\mathcal{E}.n'&=\mathcal{L}.n,\nonumber \\
\underset{q\in\{0,\ldots,Q\}}{\mathcal{E}.\nu(q)}&=
\begin{cases} \label{eq:updatenu}
\phi^{*}(i,\mathcal{L}.n,q,\mathcal{L}.\nu),&\text{if }|\mathcal{L}.\delta|=1,\\
\min\{\phi^{*}(i,\mathcal{L}.n,q,\mathcal{L}.\nu),\phi^{*}(i,\mathcal{L}.n',q,\mathcal{L}.\nu')\},&\text{if }|\mathcal{L}.\delta|>1,
\end{cases}\\
\underset{q\in\{0,\ldots,Q\}}{\mathcal{E}.\nu'(q)}&=\phi^{*}(\mathcal{L}.n,i,q,\iota(\cdot)),\label{eq:updatenup}
\end{align}
where $\oplus$ denotes sequence concatenation, and the cost-to-go function $\iota:\{0,\ldots,Q\}\mapsto\mathbb{R}$ is specified by:
\begin{equation} \nonumber
\underset{q\in\{0,\ldots,Q\}}{\iota(q)}=
\begin{cases}
c_{i0},&\text{if }|\mathcal{L}.\delta|=1,\\
\phi^{*}(i,\mathcal{L}.\mathcal{P}.n,q,\mathcal{L}.\mathcal{P}.\nu),&\text{if }|\mathcal{L}.\delta|=2,\\
\min\{\phi^{*}(i,\mathcal{L}.\mathcal{P}.n,q,\mathcal{L}.\mathcal{P}.\nu),\phi^{*}(i,\mathcal{L}.\mathcal{P}.n',q,\mathcal{L}.\mathcal{P}.\nu')\},&\text{if }|\mathcal{L}.\delta|>2.
\end{cases}
\end{equation}

When defining $\mathcal{E}.\nu$ in \eqref{eq:updatenu}, in the general case ($|\mathcal{L}.\delta|>1$) four possibilities must be considered: either customer $\mathcal{L}.n$ or $\mathcal{L}.n'$ is visited next, and in each case the visit may occur directly or after a replenishment trip. Note that the cost-to-go function changes depending on the decision. Further, when defining $\mathcal{E}.\nu'$ in \eqref{eq:updatenup}, after customer $\mathcal{L}.n$ is served, customer $i$ must be visited (directly or after replenishment). In the general case ($|\mathcal{L}.\delta|>2$), the cost-to-go function $\iota$ after $i$ is served depends on which customer, $\mathcal{L}.\mathcal{P}.n$ or $\mathcal{L}.\mathcal{P}.n'$, is visited next, and whether the visit occurs directly or after replenishment.

Finally, $\mathcal{E}.\Theta$ is computed based on the customers in $\mathcal{E}.\delta$:
\begin{equation} \nonumber
\mathcal{E}.\Theta=\kappa+\sum_{i\in\mathcal{E}.\delta}\alpha_{i}+\sum_{S\in\mathcal{J}}\left\lfloor\sum_{i\in S}[i\in\mathcal{L}.\delta]/2\right\rfloor\beta_{S}+\sum_{S\in\mathcal{C}}[S\cap\mathcal{L}.\delta\neq\emptyset]\gamma_{S}.
\end{equation}

\subsection{Completion Bounds} \label{sec:cbounds}
The labeling algorithm relies on completion bounds to eliminate labels and better control their combinatorial growth. Three distinct bounding procedures are proposed, one based on a knapsack problem, and the other two on resource-constrained shortest-path (RCSP) problems. Let $A(\mathcal{L})$ be the set of feasible routes that can be obtained by extending a label $\mathcal{L}$. More precisely, each route $r\in A(\mathcal{L})$ is such that $\mathcal{L}.\delta\oplus(0)$ is a contiguous subsequence of $r\oplus(0)$. Given a label $\mathcal{L}$, each bounding procedure defines a lower bound on the reduced cost of all routes in $A(\mathcal{L})$. When this bound is nonnegative, $\mathcal{L}$ can be discarded since it is not able to generate a profitable column.

\subsubsection{Elementary Knapsack Bound}
Consider a label $\mathcal{L}$. Route $\mathcal{L}.\delta$ has the minimum cost among all routes in $A(\mathcal{L})$, because the switch policy allows preventive restocking, and so adding customers between the depot and customer $\mathcal{L}.n$ cannot decrease costs. The knapsack bound consists in subtracting from $\overline{C}_{\mathcal{L}.\delta}$ an upper bound on the reduced cost decrease promoted by further extending $\mathcal{L}$. We denote by $\overline{\theta}_{i}$ the maximum decrease in reduced cost when extending $\mathcal{L}$ to a customer $i\notin\mathcal{L}.\delta$:
\begin{equation} \label{eqincredcost}
\overline{\theta}_{i}=\alpha_{i}+\sum_{S\in\mathcal{C}}[i\in S]\gamma_{S}.
\end{equation}

Since dual values $\beta_{S}$, $S\in\mathcal{J}$, are nonpositive, they can be ignored when computing $\overline{\theta}_{i}$. Therefore, the reduced cost of any route $r\in A(\mathcal{L})$ is bounded as follows:
\begin{equation} \nonumber
\overline{C}_{r}\geq\overline{C}_{\mathcal{L}.\delta}-\textsc{kp}(\mathcal{L}),
\end{equation}
where $\textsc{kp}(\mathcal{L})$ is the optimal solution value of a $\{0,1\}$-knapsack with set of items $\mathcal{I}=\{1,\dots,N\}\setminus\mathcal{L}.\delta$, where each item has a value of $\overline{\theta}_{i}$ and a weight of $\mathbb{E}[D_{i}]$, and the capacity of the knapsack is $fQ-\sum_{i\in\mathcal{L}.\delta}\mathbb{E}[D_{i}]$.

\textbf{Complexity.} The $\{0,1\}$-knapsack can be efficiently solved by dynamic programming in $\mathcal{O}(fQN)$ time. The knapsack bound is applied to all labels except those already pruned by RCSP-based bounds, which are evaluated first since they have a lower complexity. It is certainly computationally demanding to solve a knapsack problem at most labels. However, in Section \ref{sec:expbcnp} we verify that the overall performance of the pricing algorithm improves when knapsack bounds are active.

\subsubsection{RCSP-based Bounds}
Given a label $\mathcal{L}$, these completion bounds are computed based on non-elementary RCSPs from the depot to nodes $\mathcal{L}.n$ and $\mathcal{L}.n'$. We start with some definitions. For a given customer set $S$, let $G'(S)=(V,A')$ be an auxiliary graph with set of nodes $V$ (as in our original graph $G$) and set of arcs $A'=\{(i,j)\in A:i\notin S\}$. With each arc $(i,j)\in A'$ a modified cost $\overline{c}_{ij}$ is associated:
\begin{equation} \nonumber
\overline{c}_{ij}=
\begin{cases}
c_{ij}-\overline{\theta}_{j},&\text{if }j\neq 0,\\
c_{ij},&\text{if }j=0.\\
\end{cases}
\end{equation}

We denote by $\textsc{sp}(S,i,R)$ the value of the RCSP in $G'(S)$ from node 0 to node $i$, where $R$ resource units are initially available, and each visit to node $i\in V\setminus\{0\}$ consumes $\mathbb{E}[D_{i}]$ resource units.

Next, consider a label $\mathcal{L}$ such that $\mathcal{L}.\mathcal{P}\neq\varnothing$ and $\mathcal{L}.\mathcal{P}.\mathcal{P}\neq\varnothing$. We are interested in finding a bound on the reduced cost of all routes in $A(\mathcal{L})$. For this purpose, we let $r^{*}$ be the route with minimum reduced cost among all routes in $A(\mathcal{L})$, and define a partition $\{\{N_{k}\},k\in\{-1,0,1\}\}$ of the sample space, in which $N_{0}$ is the event where customer $\mathcal{L}.n$ is visited in its actual position in $r^{*}$, and $N_{-1}$ and $N_{1}$ are the events where $\mathcal{L}.n$ is visited in the immediately preceding and following positions, respectively. The reduced cost of $r^{*}$ can be rewritten and bounded as follows:
\begin{align}
\overline{C}_{r^{*}}&=\mathbb{E}[C_{r^{*}}]-\Theta_{r^{*}}, \nonumber \\
&=\mathbb{E}[C_{r^{*}}-\Theta_{r^{*}}], \nonumber \\
&\geq\min_{k\in\{1,0,-1\}}\mathbb{E}[C_{r^{*}}-\Theta_{r^{*}}\mid N_{k}], \label{eqcondexp}
\end{align}
where $-\Theta_{r^{*}}$ is the contribution of the dual values to the reduced cost of $r^{*}$ (see equation \eqref{eqredcost}).

We now derive lower bounds on the conditional expectations in \eqref{eqcondexp}, starting with event $N_{0}$. When $\mathcal{L}.n$ is visited in its planned position, a lower bound on the remaining cost once $\mathcal{L}.n$ is served is given by $\mathcal{L}.\nu(Q)$. The contribution of all dual values associated with $\mathcal{L}.\delta$ to the reduced cost is given by $-\mathcal{L}.\Theta$. Furthermore, the reduced cost change promoted by extending $\mathcal{L}$ is bounded from below by the value of the RCSP in $G'(\mathcal{L}.\delta)$ from the depot to customer $\mathcal{L}.n$ subject to a resource limit of $fQ-\sum_{i\in\mathcal{L}.\delta\setminus\{\mathcal{L}.n\}}\mathbb{E}[D_{i}]$. Hence, the following inequality holds:
\begin{equation} \label{eqcondn0}
\mathbb{E}[C_{r^{*}}-\Theta_{r^{*}}\mid N_{0}]\geq\textsc{sp}\Big(\mathcal{L}.\delta,\mathcal{L}.n,fQ-\!\!\!\!\!\!\!\sum_{i\in\mathcal{L}.\delta\setminus\{\mathcal{L}.n\}}\!\!\!\!\!\!\!\mathbb{E}[D_{i}]\Big)+\overline{\theta}_{\mathcal{L}.n}+\mathcal{L}.\nu(Q)-\mathcal{L}.\Theta,
\end{equation}
where the term $\overline{\theta}_{\mathcal{L}.n}$ is added back because the contribution of the duals associated with $\mathcal{L}.n$ is already taken into account in $\mathcal{L}.\Theta$.

Similarly, event $N_{1}$ corresponds to the case where customers $\mathcal{L}.n$ and $\mathcal{L}.n'$ are switched in the a posteriori route. In this case, $\mathcal{L}.\nu'(Q)$ is a lower bound on the remaining cost once $\mathcal{L}.n'$ is served, and a bound on the reduced cost of $r^{*}$ is given by:
\begin{equation} \label{eqcondn1}
\mathbb{E}[C_{r^{*}}-\Theta_{r^{*}}\mid N_{1}]\geq\textsc{sp}\Big(\mathcal{L}.\delta,\mathcal{L}.n',fQ-\!\!\!\!\!\!\!\sum_{i\in\mathcal{L}.\delta\setminus\{\mathcal{L}.n'\}}\!\!\!\!\!\!\!\mathbb{E}[D_{i}]\Big)+\overline{\theta}_{\mathcal{L}.n'}+\mathcal{L}.\nu'(Q)-\mathcal{L}.\Theta.
\end{equation}

The last case, which corresponds to event $N_{-1}$, is when customer $\mathcal{L}.n$ is visited in the previous position. In this case, the remaining cost once $\mathcal{L}.n$ is served depends on the customer $j$ with which $\mathcal{L}.n$ is switched. A lower bound $l(\mathcal{L},j)$ on such cost is given by:
\begin{equation} \label{eqlbLj}
l(\mathcal{L},j)=c_{\mathcal{L}.n,j}+\min\{c_{j,\mathcal{L}.\mathcal{P}.n}+\mathcal{L}.\mathcal{P}.\nu(Q),c_{j,\mathcal{L}.\mathcal{P}.n'}+\mathcal{L}.\mathcal{P}.\nu'(Q)\}.
\end{equation}

The first term in \eqref{eqlbLj} assumes that $j$ is visited directly after $\mathcal{L}.n$, and the second term is a lower bound on the cost-to-go once $j$ is served. Then, all unvisited customers $j\notin\mathcal{L}.\delta$ must be considered when computing the bound on the expectation conditional on $N_{-1}$:
\begin{equation} \label{eqcondnm1}
\begin{split}
\mathbb{E}[C_{r^{*}}-\Theta_{r^{*}}\mid N_{-1}]\geq\min_{j\notin\mathcal{L}.\delta,j\neq 0}\bigg\{\textsc{sp}\Big(\mathcal{L}.\delta\cup\{j\},\mathcal{L}.n,fQ\\-\!\!\!\!\!\!\!\sum_{i\in(\mathcal{L}.\delta\cup\{j\})\setminus\{\mathcal{L}.n\}}\!\!\!\!\!\!\!\mathbb{E}[D_{i}]\Big)+l(\mathcal{L},j)-\overline{\theta}_{j}\bigg\}+\overline{\theta}_{\mathcal{L}.n}-\mathcal{L}.\Theta,
\end{split}
\end{equation}
where $\overline{\theta}_{j}$ is explicitly subtracted because it is not taken into account in $\mathcal{L}.\Theta$.

Therefore, by \eqref{eqcondexp} the minimum among \eqref{eqcondn0}, \eqref{eqcondn1} and \eqref{eqcondnm1} is a lower bound on the reduced cost of all routes in $A(\mathcal{L})$, that is, a completion bound.

\textbf{Bound variant.} The bound derived above loses effectiveness as more and more strengthened capacity cuts are separated. The reason is that the dual value of a capacity cut may be counted multiple times in the computation of RCSPs. To alleviate this issue, we propose a variant of the bound in which dual values of capacity cuts are disregarded when computing \eqref{eqcondn0}, \eqref{eqcondn1} and \eqref{eqcondnm1}, and the total contribution of those cuts, given by $\sum_{S\in\mathcal{C}}\gamma_{S}$, is subtracted after computing the minimum among \eqref{eqcondn0}, \eqref{eqcondn1} and \eqref{eqcondnm1}. This assumes a worst-case scenario where all capacity cuts are active, but avoids multiple counting of capacity cut duals.

\textbf{Complexity.} In \eqref{eqcondn0}, \eqref{eqcondn1} and \eqref{eqcondnm1}, the (non-elementary) RCSPs can be obtained by dynamic programming in $\mathcal{O}(fQN^2)$ time. Because of the minimum operation in \eqref{eqcondnm1}, the runtime complexity for computing RCSP-based bounds is $\mathcal{O}(fQN^3)$. This is clearly too demanding to apply to every label generated by the pricing algorithm. For this reason, we employ the speed-up procedure proposed in \cite{florioevrp2021}. Before starting the labeling algorithm, we select a subset of customers $M\subset\{1,\ldots,N\}$ and pre-compute $\textsc{sp}(S,i,R)$ for all $S\subset M$, $i\in\{1,\ldots,N\}$, and $R\in\{0,\ldots,fQ\}$. Then, whenever a RCSP is required for computing bounds, say $\textsc{sp}(S,i,R)$, we retrieve the pre-computed value $\textsc{sp}(S\cap M,i,R)$ (in $\mathcal{O}(1)$ time). As a result, bound evaluation occurs in $\mathcal{O}(N)$ time, at the expense of less tight bounds. The pre-computation stage adds a runtime complexity of $\mathcal{O}(2^{|M|}f^2Q^2N^3)$ each time the pricing algorithm is invoked. In our implementation, as in \cite{florioevrp2021}, we consider in the customer set $M$, at each pricing iteration, the eight customers with largest $\alpha_{i}/\mathbb{E}[D_{i}]$ ratios.

\subsection{Cut Separation and Branching Rules}
Violated capacity inequalities are identified using the CVRPSEP package by \cite{lysgaard2004new}, and added to the RMP as strengthened capacity cuts. Subset row cuts defined over customer triplets are separated by enumeration. The RMP is solved by alternating between column and row generation until an optimal solution is found. In traditional branch-cut-and-price, non-robust cuts must be separated cautiously to avoid pricing intractability \citep{poggi2014new}. In our case, however, this is not an issue because we do not employ dominance when solving the pricing problem. When the solution to the RMP is fractional, branching is performed on the variables of the arc-based formulation. We apply a semi-strong branching rule, in which all variables are considered, and branching is performed on the variable that leads to the highest increase in the lower bound, when considering only the columns and cuts generated so far. Finally, in order to obtain good upper bounds during the execution of the algorithm, the set-partitioning formulation is solved exactly with the columns generated thus far whenever a node of the branch-and-bound tree is fully processed.

\section{Computational Results} \label{sec:results}
Two sets of computational experiments are performed. In Section \ref{sec:exppolcomp}, we assume a planned route is given and evaluate the detour-to-depot, optimal restocking and switch policies under different demand distributions and load factors. Next, in Section \ref{sec:expbcnp}, we apply the branch-cut-and-price algorithm from Section \ref{sec:bcnpalg} to solve VRPSD instances under the switch policy, and compare results with those obtained by solving the same instances under optimal restocking. The algorithm is implemented in C++, and IBM\textsuperscript{\textregistered} CPLEX\textsuperscript{\textregistered} version 12.10 is employed to solve linear programs. All experiments were performed on a single core of an Intel\textsuperscript{\textregistered} Xeon\textsuperscript{\textregistered} Gold 6130 (2.1GHz) processor with 12GB of available RAM.

\subsection{Policy Comparison for Non-optimal Planned Sequences} \label{sec:exppolcomp}
Solving the VRPSD to optimality is a considerable challenge. In practice, one may reasonably expect to obtain a good VRPSD solution by first computing a set of planned routes with a deterministic VRP algorithm, and then executing these routes under the desired recourse policy. As discussed in \cite{florio2020optimal}, this strategy leads to good solutions to the single-vehicle version of the problem. The main goal of this section is to compare different recourse policies when planned sequences are not optimized for any recourse policy. In particular, we assess the relative benefits of the switch policy over optimal restocking under a range of demand distributions, load factors, sequence sizes (i.e., number of customers) and penalty values for failures.

The instances used in these experiments are randomly generated. Customer locations are uniformly distributed on a 1,000 by 1,000 grid, and the depot is located at the corner $(0,0)$. This network layout is the same as in \cite{SecomandiMargot2009}, with the depot located in the south-west corner. Three levels of demand variability are considered. In the low variability level, demands follow binomial distributions with a variance-to-mean ratio of 0.5. In the medium level, demands are Poisson distributed (variance-to-mean ratio of 1). In the high variability level, demands follow negative binomial distributions with a variance-to-mean ratio of 2. In all distributions, left and right tail probabilities below $10^{-6}$ are truncated to zero. The expected demand of each customer is randomly drawn from a discrete uniform distribution over the interval $[10,100]$. The chosen distributions allow us to experiment with different demand variabilities while keeping the same mean demands. Moreover, considering also the range from which mean demand values are drawn, these distributions feature extensive support sets and thus lead to highly stochastic instances.

A planned route $r$ is obtained by solving a traveling salesman problem on the set of generated customers and the depot. The number of customers in each route $r$, indicated by $|r|$, ranges from $3$ to $15$. Concerning load factors, we experiment with values $f\in\{1.3,1.6,1.9,2.5\}$ by setting, for each route $r$, the vehicle capacity to $\sum_{i\in r}\mathbb{E}[D_{i}]/f$ rounded to the nearest integer. With these load factors, one or two restocking trips are required in most cases to fully serve all demand. We generate $20$ instances for each variability level and values of $|r|$ and $f$, leading to a total of $20\times3\times13\times4=3,120$ instances.

\begin{table}[t]
\centering
\begingroup
\small
\renewcommand{\arraystretch}{0.825}
\caption{\label{tab:policycmp}Switch Policy: Improvements over Optimal Restocking (Poisson Demands)}
\begin{tabular}{rrrrrrrrrr}
\toprule
& \multicolumn{2}{c}{$f=1.3$} & \multicolumn{2}{c}{$f=1.6$} & \multicolumn{2}{c}{$f=1.9$}  & \multicolumn{2}{c}{$f=2.5$} \\
\cmidrule(l){2-3}\cmidrule(l){4-5}\cmidrule(l){6-7}\cmidrule(l){8-9}
$|r|$ & Avg(\%) & Max(\%) & Avg(\%) & Max(\%) & Avg(\%) & Max(\%) & Avg(\%) & Max(\%) & Avg \\
\midrule
3 & 0.19 & 2.45 & 1.26 & 4.51 & 1.50 & 6.40 & 0.39 & 1.42 & 0.83 \\
4 & 0.27 & 2.72 & 1.43 & 13.63 & 1.73 & 9.18 & 1.12 & 7.30 & 1.14 \\
5 & 0.06 & 1.24 & 0.26 & 1.67 & 0.81 & 3.61 & 0.57 & 3.10 & 0.43 \\
6 & 0.57 & 3.22 & 1.35 & 7.34 & 2.56 & 10.23 & 1.10 & 4.55 & 1.40 \\
7 & 0.13 & 1.40 & 0.60 & 4.15 & 1.64 & 7.25 & 0.73 & 4.36 & 0.77 \\
8 & 0.57 & 3.05 & 0.52 & 3.12 & 0.56 & 5.86 & 0.88 & 4.80 & 0.63 \\
9 & 0.39 & 2.73 & 1.07 & 5.79 & 1.07 & 4.55 & 0.81 & 3.08 & 0.84 \\
10 & 0.33 & 2.02 & 0.27 & 1.48 & 0.90 & 7.28 & 0.55 & 1.70 & 0.51 \\
11 & 0.15 & 1.42 & 0.59 & 3.74 & 0.63 & 4.20 & 0.58 & 2.41 & 0.49 \\
12 & 0.38 & 1.61 & 0.38 & 1.30 & 0.66 & 3.18 & 0.69 & 2.38 & 0.53 \\
13 & 0.42 & 3.75 & 0.50 & 2.20 & 0.56 & 3.60 & 0.55 & 2.44 & 0.51 \\
14 & 0.42 & 3.64 & 0.71 & 4.74 & 1.10 & 4.99 & 0.60 & 2.70 & 0.71 \\
15 & 0.17 & 0.94 & 0.81 & 3.23 & 0.90 & 4.75 & 0.86 & 3.51 & 0.68 \\
Avg & 0.31 &  & 0.75 &  & 1.12 &  & 0.72 &  & 0.73 \\
\bottomrule
\end{tabular}
\endgroup
\end{table}

Table \ref{tab:policycmp} lists the results obtained by considering Poisson distributions. Columns ``Avg(\%)'' and ``Max(\%)'' indicate the average and maximum percentage improvements on expected cost, respectively, when applying the switch policy instead of the optimal restocking policy. These maximum and average values are taken over all $20$ instances for each pair of $|r|$ and $f$ values. Note that the expected cost of a route depends on its orientation. Cost improvements are calculated assuming the route orientation that minimizes the expected cost, for each policy. Savings enabled by the switch policy are, on average, 0.73\% over all 1,040 instances. In several cases, however, improvements exceed 5\%, and values above 10\% are observed in a couple of instances. These improvements are not significantly correlated with the number of customers nor with the load factor, although gains tend to be smaller when $f=1.3$, the lowest load factor experimented.

\begin{figure}[t]
\centering
	{\includegraphics[trim={0 90 460 0},clip,width=300pt]{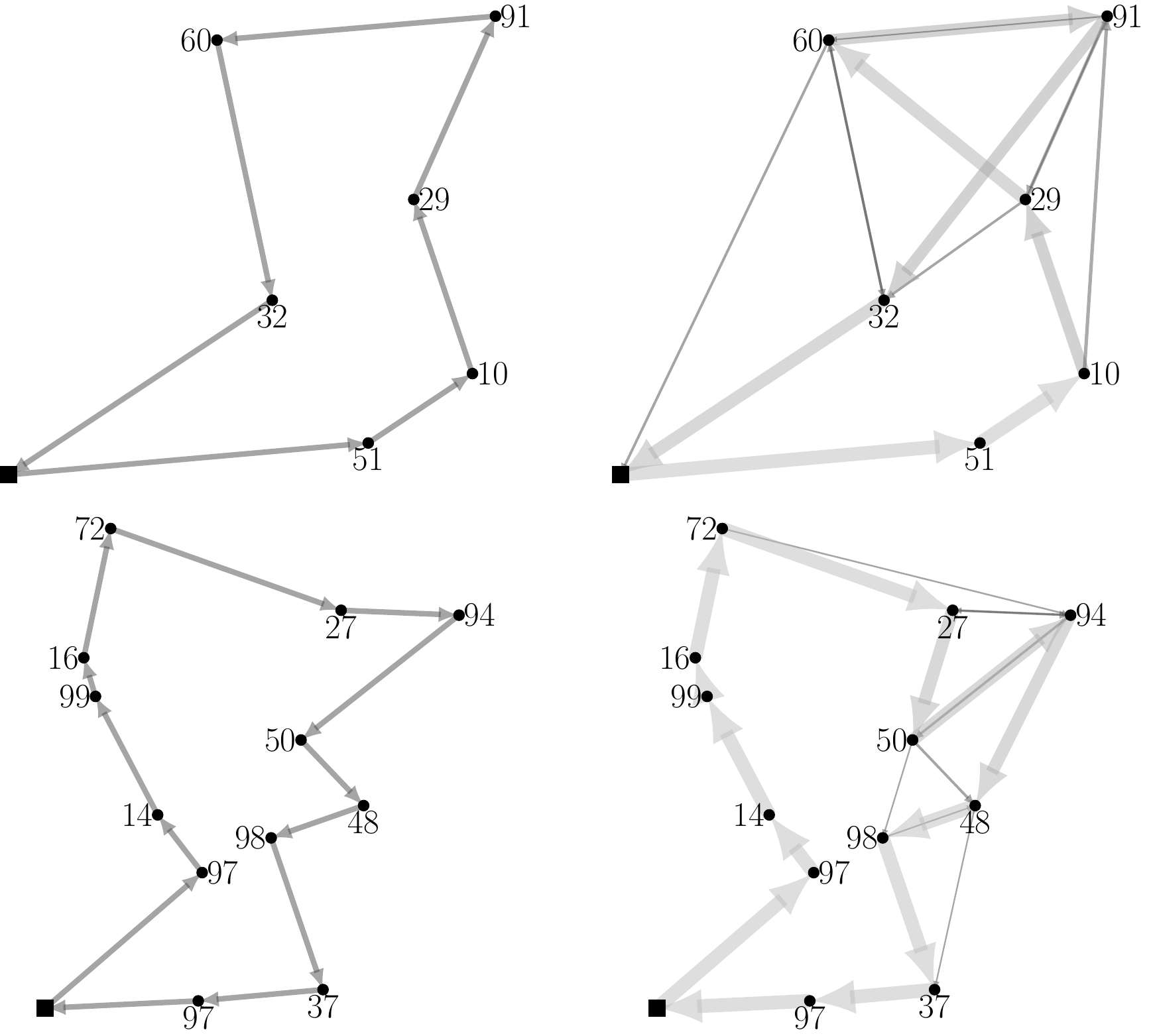}}
	{}
		\caption{Planned sequences (left drawings) and simulated a posteriori sequences (right drawings) under the switch policy. Improvements upon optimal restocking are $7.3\%$ and $3.2\%$ in the top and bottom instances, respectively.\label{fig:policycmp}}
\end{figure}

Figure \ref{fig:policycmp} illustrates two planned routes (left-hand side drawings) and their corresponding a posteriori executions under the switch policy (right-hand side drawings) obtained by simulation over 5,000 scenarios. The values next to the customers indicate their expected demand. In the right drawings, the thickness of the arrows is proportional to the probability of the connected nodes being visited consecutively, and restocking trips are omitted for simplification. In the top instance, the vehicle capacity is $Q=170$ and the load factor is $f\approx1.6$. In the bottom instance, $Q=394$ and $f\approx1.9$. Savings over optimal restocking are $7.3\%$ and $3.2\%$ in the top and bottom instances, respectively. Note that the most probable visiting sequences (indicated by the paths formed by the thickest arrows) deviate from the planned sequences, and that customers are frequently switched to reduce expected costs.

Next, in Table \ref{tab:policycmp2}, we compare the detour-to-depot and switch policies with optimal restocking for varying load factors and demand distributions. As observed, savings are larger when moving from detour-to-depot to optimal restocking, when compared to moving from optimal restocking to the switch policy. Interestingly, increasing demand variability is not associated with larger improvements when moving to a more sophisticated policy. In fact, in most cases savings are higher when demand variability is low, which, as we interpret, enables more sensible restocking and customer switching decisions.

\begin{table}[t]
\centering
\begingroup
\small
\renewcommand{\arraystretch}{0.825}
\caption{\label{tab:policycmp2}Policy Comparison for Varying Load Factors and Demand Variability Levels}
\begin{tabular}{rrrrrrrrr}
\toprule
& \multicolumn{4}{c}{Detour-Depot \textit{vs} Opt. R.(\%)$^a$} & \multicolumn{4}{c}{Switch \textit{vs} Opt. R.(\%)} \\
\cmidrule(l){2-5}\cmidrule(l){6-9}
$f$ & L$^b$ & M$^b$ & H$^b$ & Avg & L & M & H & Avg \\
\midrule
1.3 & $-$3.10 & $-$3.05 & $-$2.94 & $-$3.03 & 0.33 & 0.31 & 0.31 & 0.31 \\
1.6 & $-$3.95 & $-$3.79 & $-$3.54 & $-$3.76 & 0.76 & 0.75 & 0.74 & 0.75 \\
1.9 & $-$2.05 & $-$2.32 & $-$2.64 & $-$2.34 & 1.26 & 1.12 & 0.98 & 1.12 \\
2.5 & $-$4.80 & $-$4.61 & $-$4.28 & $-$4.56 & 0.75 & 0.72 & 0.67 & 0.72 \\
Avg & $-$3.47 & $-$3.44 & $-$3.35 & $-$3.42 & 0.77 & 0.73 & 0.67 & 0.72 \\
\bottomrule
\end{tabular}
\captionsetup{justification=centering}
\caption*{$^a$ Negative values indicate worse performance (higher costs) than optimal restocking. \\$^b$ `L', `M', `H' denote low, medium and high levels of demand variability.}
\endgroup
\end{table}

\begin{figure}[b]
\centering
	{\includegraphics[trim={0 0 0 0},clip,width=\textwidth]{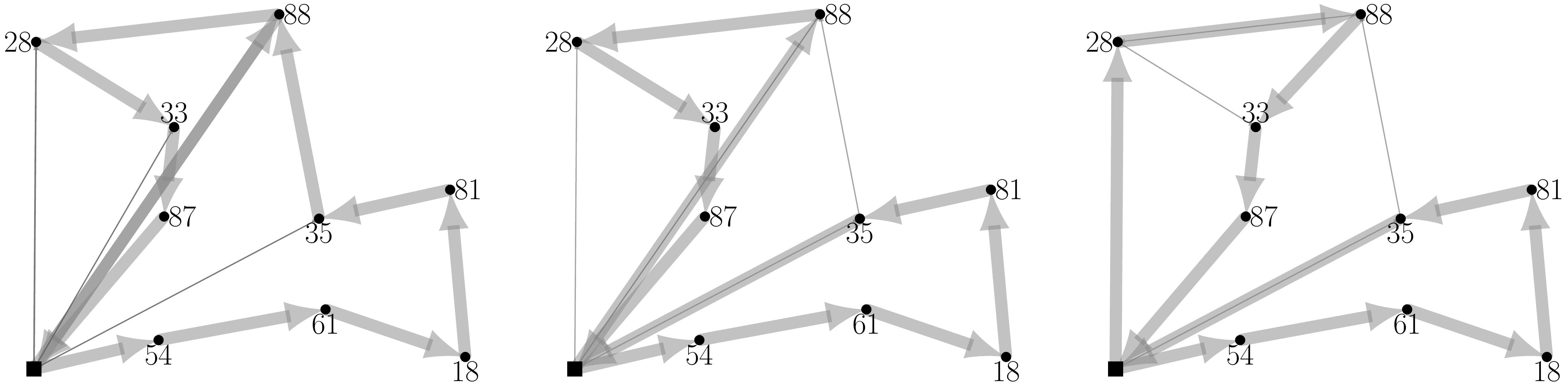}}
	{}
		\caption{Simulated a posteriori sequences under the detour-to-depot (left), optimal restocking (center) and switch (right) policies. In this instance, the switch policy enables savings of $5.8\%$ over optimal restocking.\label{fig:policycmp2}}
\end{figure}

Figure \ref{fig:policycmp2} compares the simulated a posteriori sequences produced by the detour-to-depot, optimal restocking and switch policies on the same planned route. Replenishment trips to the depot are also included in the drawings. The vehicle capacity is $Q=303$ and the load factor is $f\approx1.6$. Under detour-to-depot, with high probability a failure occurs when visiting the sixth customer, which has an expected demand of 88 (with lower probability, as indicated by the thin arrows, a failure occurs at the preceding or at the following two customers instead). Under optimal restocking, a preventive replenishment trip is usually performed after visiting the fifth customer. Finally, under the switch policy, a preventive replenishment trip is performed \emph{and} the sixth and seventh customers are switched. In this instance, under Poisson demands, the switch policy enables savings of $5.8\%$ over optimal restocking, which, in turn, is $2.7\%$ less costly than the detour-to-depot policy.

\begin{figure}[t]
\centering
	{\includegraphics[trim={0 10 0 0},clip,width=270pt]{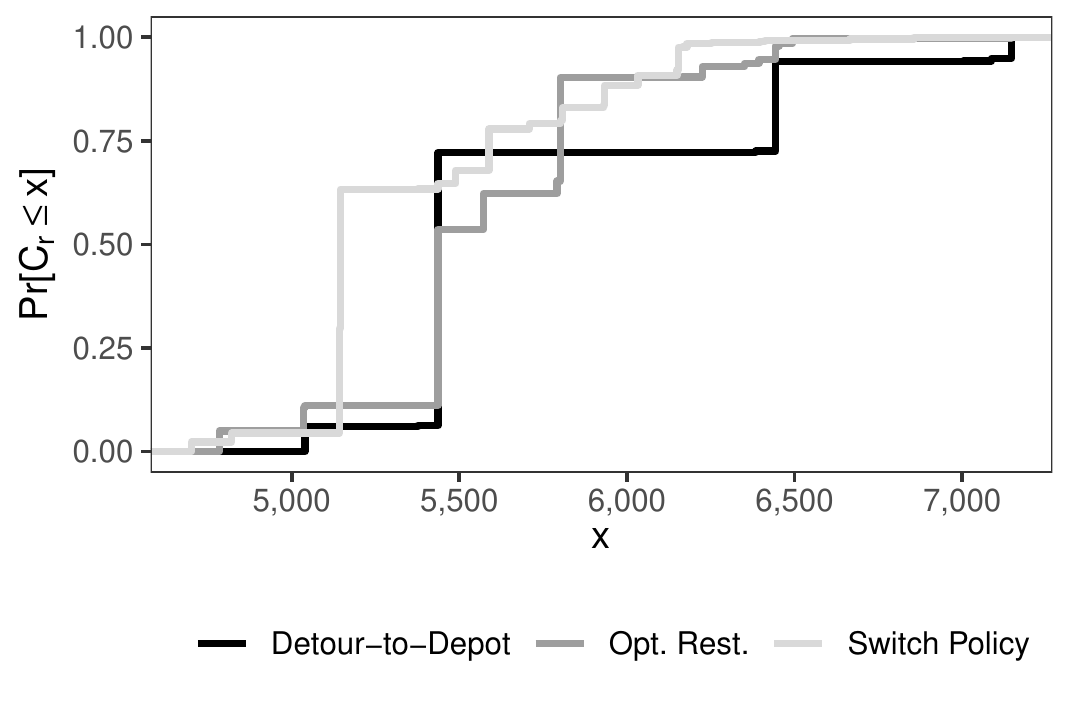}}
	{}
		\caption{Empirical cumulative distribution functions of the cost of a planned route under the detour-to-depot, optimal restocking and switch policies. The augmented decision space of the switch policy leads to many more potential a posteriori sequences.\label{fig:policycmp3}}
\end{figure}

To further develop insights into the differences among policies, we analyze the empirical cumulative distribution functions (CDFs) of the cost of a route under the three policies. As an example, we consider a planned route $r$ with $6$ customers, load factor $f\approx1.9$ and high demand variability. The CDFs are depicted in Figure \ref{fig:policycmp3}. The area above the curve of each policy corresponds to the expected cost of $r$ under that policy. The higher number of steps in the CDF associated with the switch policy is a consequence of the much larger decision space of that policy, which leads to many more potential a posteriori sequences. Note that the CDF of the simpler detour-to-depot policy consists of far fewer steps.

As mentioned in Remark \ref{remarkpenalty}, in some cases it may be desirable to discourage split deliveries by penalizing failures. Recall that a failure occurs when arriving at a customer with insufficient load to serve that customer. We define failure rate as the expected number of failures when executing a planned route, and denote by $\Lambda_{\mathtt{OR}}$ and $\Lambda_{\mathtt{SP}}$ the failure rates under the optimal restocking and switch policies, respectively. In Table \ref{tab:policycmp3}, we compare the failure rates of both policies when a penalty $p$ is applied for each failure. We also report the percentage savings of the switch policy over optimal restocking in column ``Impr.(\%)''. Again, the switch policy demonstrates a slight edge against optimal restocking, both in terms of expected cost and rate of failures.

\vspace{6pt}
\begin{table}[h]
\centering
\begingroup
\small
\renewcommand{\arraystretch}{0.825}
\caption{\label{tab:policycmp3}Failure Rates of the Optimal Restocking and Switch Policies for Varying Penalty Values}
\begin{tabular}{rrrrrrrrrr}
\toprule
& \multicolumn{3}{c}{$p=10$} & \multicolumn{3}{c}{$p=100$} & \multicolumn{3}{c}{$p=1,000$} \\
\cmidrule(l){2-4}\cmidrule(l){5-7}\cmidrule(l){8-10}
$f$ & Impr.(\%) & $\Lambda_{\mathtt{OR}}$ & $\Lambda_{\mathtt{SP}}$ & Impr.(\%) & $\Lambda_{\mathtt{OR}}$ & $\Lambda_{\mathtt{SP}}$ & Impr.(\%) & $\Lambda_{\mathtt{OR}}$ & $\Lambda_{\mathtt{SP}}$ \\
\midrule
1.3 & 0.32 & 0.105 & 0.094 & 0.33 & 0.065 & 0.058 & 0.42 & 0.018 & 0.015 \\
1.6 & 0.75 & 0.171 & 0.150 & 0.79 & 0.127 & 0.104 & 1.15 & 0.058 & 0.043 \\
1.9 & 1.11 & 0.626 & 0.661 & 1.04 & 0.501 & 0.522 & 1.13 & 0.117 & 0.106 \\
2.5 & 0.72 & 0.681 & 0.647 & 0.79 & 0.543 & 0.481 & 1.39 & 0.237 & 0.205 \\
Avg & 0.72 & 0.396 & 0.388 & 0.74 & 0.309 & 0.292 & 1.02 & 0.108 & 0.092 \\
\bottomrule
\end{tabular}
\captionsetup{justification=centering}
\caption*{Note: the values in columns $\Lambda_{\mathtt{OR}}$ and $\Lambda_{\mathtt{SP}}$ correspond to average failure rates over all $3,120$ instances,\\each simulated over $5,000$ scenarios.}
\endgroup
\end{table}

\subsection{Branch-Cut-and-Price} \label{sec:expbcnp}
We apply the branch-cut-and-price algorithm on 21 instances with up to 50 nodes of the sets $\mathtt{A}$ and $\mathtt{P}$ of the CVRPLIB \citep{uchoa2017new}. Each instance is solved under three load factors $f\in\{1.3,1.6,1.9\}$ and the same demand distributions defined in Section \ref{sec:exppolcomp}. Hence, we run the exact algorithm on a total of $21\times3\times3=189$ instances. In each instance, the number of available vehicles is set to the minimum required for feasibility, and the vehicle capacity is adjusted to $Q'/q_f$ (rounded to the nearest integer), where $Q'$ is the original vehicle capacity, $q_{1.3}=2$, $q_{1.6}=3$, and $q_{1.9}=4$. This capacity adjustment compensates the increase in the average number of customers served per route caused by load factors larger than 1. It also enables our algorithm to find optimal solutions under large load factors, as instance difficulty increases with $f$ even when $fQ$ remains constant. We set a maximum runtime limit of 2 hours for each instance.

Table \ref{tab:summbcnp} presents a summary of the results obtained by the algorithm. Each group of columns corresponds to a load factor, and the results are aggregated over the three demand variability settings. Columns ``LB'', ``Opt'', ``Gap'' and ``Time'' indicate the number of instances where the lower bound at the root node could be computed, the number of instances solved to optimality, the average percentage optimality gap in the unsolved instances with known lower bounds, and the average runtime in minutes, respectively. Out of 189 instances, 90 could be solved to optimality. For 35 instances, however, the linear relaxation at the root node of the branch-and-bound tree could not be computed. In these instances, the completion bounds are not effective enough to control the combinatorial growth of labels when solving the pricing problem. The average optimality gap in the remaining 64 instances is only 0.34\%. We attribute this low gap to elementary pricing together with an aggressive separation of strengthened capacity and subset row cuts. The detailed results for each instance are reported in Tables \ref{tab:bcpall13}-\ref{tab:bcpall19}, in the appendix.

\begin{table}[t]
\centering
\begingroup
\small
\renewcommand{\arraystretch}{0.825}
\caption{\label{tab:summbcnp}Summarized Branch-Cut-and-Price Results over 189 Instances}
\begin{tabular}{lrrrrrrrrrrrrrr}
\toprule
& \multicolumn{4}{c}{$f=1.3$} &&\multicolumn{4}{c}{$f=1.6$} && \multicolumn{4}{c}{$f=1.9$} \\
\cmidrule(l){2-5}\cmidrule(l){7-10}\cmidrule(l){12-15}
Instance & LB & Opt & Gap & Time && LB & Opt & Gap & Time && LB & Opt & Gap & Time \\
\midrule
A-n32-k5 & 3 & 2 & 1.17 & 86.4 && 3 & 0 & 0.23 & 120.0 && 3 & 3 & --- & 47.7 \\
A-n33-k5 & 3 & 3 & --- & 17.8 && 3 & 2 & 0.00 & 70.8 && 3 & 3 & --- & 16.0 \\
A-n34-k5 & 3 & 1 & 0.11 & 104.6 && 3 & 3 & --- & 41.1 && 3 & 3 & --- & 22.7 \\
A-n36-k5 & 0 & 0 & --- & 120.0 && 3 & 0 & 0.61 & 120.0 && 3 & 2 & 0.02 & 52.7 \\
A-n37-k5 & 1 & 1 & --- & 104.6 && 1 & 0 & 0.34 & 120.0 && 0 & 0 & --- & 120.0 \\
A-n38-k5 & 3 & 1 & 0.19 & 90.4 && 3 & 2 & 0.06 & 77.4 && 3 & 1 & 0.13 & 83.3 \\
A-n39-k5 & 0 & 0 & --- & 120.0 && 3 & 0 & 0.48 & 120.0 && 3 & 2 & 0.14 & 85.0 \\
A-n44-k6 & 1 & 0 & 1.36 & 120.0 & &3 & 0 & 0.40 & 120.0 && 3 & 2 & 0.08 & 74.9 \\
A-n45-k6 & 2 & 1 & 0.01 & 101.1 && 3 & 0 & 0.16 & 120.0 && 3 & 3 & --- & 11.7 \\
A-n45-k7 & 0 & 0 & --- & 120.0 && 3 & 0 & 0.18 & 120.0 && 3 & 0 & 0.02 & 120.0 \\
A-n46-k7 & 0 & 0 & --- & 120.0 && 2 & 0 & 0.82 & 120.0 && 3 & 1 & 0.15 & 91.3 \\
A-n48-k7 & 0 & 0 & --- & 120.0 && 1 & 0 & 0.27 & 120.0 && 3 & 1 & 0.11 & 91.1 \\
P-n19-k2 & 3 & 3 & --- & 4.0 && 3 & 2 & 0.22 & 54.6 && 3 & 3 & --- & 2.6 \\
P-n20-k2 & 3 & 3 & --- & 12.7 && 3 & 1 & 0.71 & 113.9 && 3 & 3 & --- & 5.4 \\
P-n21-k2 & 3 & 2 & 0.15 & 70.6 && 3 & 3 & --- & 10.6 && 3 & 3 & --- & 2.6 \\
P-n22-k2 & 3 & 3 & --- & 48.2 && 3 & 0 & 0.60 & 120.0 && 3 & 3 & --- & 12.6 \\
P-n40-k5 & 1 & 0 & 0.76 & 120.0 && 3 & 0 & 0.35 & 120.0 && 3 & 1 & 0.11 & 89.7 \\
P-n45-k5 & 0 & 0 & --- & 120.0 && 1 & 0 & 0.79 & 120.0 && 3 & 3 & --- & 77.0 \\
P-n50-k10 & 3 & 1 & 0.26 & 88.4 && 3 & 3 & --- & 2.0 && 3 & 3 & --- & 0.1 \\
P-n50-k7 & 3 & 0 & 0.49 & 120.0 && 3 & 2 & 0.38 & 84.6 && 3 & 2 & 0.05 & 59.4 \\
P-n50-k8 & 3 & 3 & --- & 4.2 && 3 & 3 & --- & 5.0 && 3 & 3 & --- & 2.8 \\
Total & 38 & 24 & 0.43 & 86.3 && 56 & 21 & 0.41 & 90.5 && 60 & 45 & 0.09 & 50.9 \\
\bottomrule
\end{tabular}
\captionsetup{justification=centering}
\caption*{Note: columns ``Gap'' indicate the average percentage optimality gap in the unsolved instances where a\\lower bound was computed.}
\endgroup
\end{table}

In Table \ref{tab:cmpOR}, we compare the optimal solutions obtained under the switch policy and optimal restocking, considering medium demand variability (Poisson distributed demands). The solutions under optimal restocking are computed using the algorithm by \cite{Florio_2020}. Improvements enabled by the switch policy are smaller in these instances, with a maximum of $0.47\%$ and an average of $0.1\%$. The routes in the optimal solution under the switch policy are mostly different from the optimal restocking solution routes, as indicated in column ``P.d.(\%)''. We also compare the optimal solutions under both policies with respect to a ``late service'' measure defined as the number of customers served after the vehicle replenishes. In practice, replenishment may take a significant amount of time and delay service at customers visited after restocking. Column ``L.s.(\%)'' reports the reduction in the number of customers served after restocking, which is calculated by simulating each solution route over 5,000 scenarios. As shown, the switch policy also performs better on this measure and reduces late service by 1.8\%, on average.

\begin{table}[t]
\centering
\begingroup
\small
\renewcommand{\arraystretch}{0.825}
\caption{\label{tab:cmpOR}Comparison of Optimal VRPSD Solutions under Optimal Restocking and Switch Policies}
\begin{tabular}{lrrrrrrrrr}
\toprule
&  &  & \multicolumn{2}{c}{Opt. Restocking}  & \multicolumn{5}{c}{Switch Policy} \\
\cmidrule(l){4-5}\cmidrule(l){6-10}
Instance & $Q$ & $f$ & Best & Time & Best & Time & Impr.(\%)$^a$ & P.d.(\%)$^b$ & L.s.(\%)$^c$ \\
\midrule
A-n32-k5 & 50 & 1.3 & 1538.652 & 3.2 & 1538.008 & 80.4 & 0.04 & 86 & 0.63 \\
A-n32-k5 & 25 & 1.9 & 2618.126 & 0.5 & 2615.913 & 39.5 & 0.08 & 67 & $-$0.61 \\
A-n33-k5 & 50 & 1.3 & 1223.976 & 0.4 & 1223.964 & 30.4 & 0.00 & 86 & $-$0.01 \\
A-n33-k5 & 25 & 1.9 & 2020.024 & 0.5 & 2018.620 & 22.4 & 0.07 & 50 & 0.21 \\
A-n34-k5 & 33 & 1.6 & 1961.158 & 0.7 & 1961.076 & 59.3 & 0.00 & 78 & 0.05 \\
A-n34-k5 & 25 & 1.9 & 2430.036 & 0.3 & 2429.258 & 46.7 & 0.03 & 50 & $-$0.28 \\
A-n39-k5 & 25 & 1.9 & 2548.617 & 1.8 & 2548.391 & 42.6 & 0.01 & 73 & 5.28 \\
A-n44-k6 & 25 & 1.9 & 2921.987 & 0.7 & 2917.078 & 49.1 & 0.17 & 77 & $-$2.72 \\
A-n45-k6 & 25 & 1.9 & 2999.063 & 0.6 & 2997.605 & 8.9 & 0.05 & 39 & 2.14 \\
P-n19-k2 & 80 & 1.3 & 370.006 & 0.3 & 370.006 & 3.6 & 0.00 & 100 & 0.00 \\
P-n19-k2 & 53 & 1.6 & 469.549 & 0.5 & 469.539 & 33.4 & 0.00 & 75 & $-$0.01 \\
P-n19-k2 & 40 & 1.9 & 558.562 & 0.1 & 557.394 & 1.6 & 0.21 & 20 & 7.77 \\
P-n20-k2 & 80 & 1.3 & 375.401 & 0.6 & 375.066 & 9.4 & 0.09 & 100 & $-$0.35 \\
P-n20-k2 & 40 & 1.9 & 567.404 & 0.1 & 566.236 & 4.5 & 0.21 & 40 & 8.05 \\
P-n21-k2 & 80 & 1.3 & 363.021 & 3.0 & 363.013 & 42.6 & 0.00 & 67 & 0.47 \\
P-n21-k2 & 53 & 1.6 & 467.072 & 0.3 & 466.729 & 8.8 & 0.07 & 50 & $-$0.26 \\
P-n21-k2 & 40 & 1.9 & 560.199 & 0.1 & 558.616 & 2.9 & 0.28 & 75 & 6.26 \\
P-n22-k2 & 80 & 1.3 & 375.565 & 3.1 & 375.197 & 43.1 & 0.10 & 100 & $-$0.51 \\
P-n22-k2 & 40 & 1.9 & 576.611 & 0.2 & 573.906 & 5.0 & 0.47 & 100 & 15.61 \\
P-n45-k5 & 38 & 1.9 & 1269.893 & 2.4 & 1266.367 & 103.9 & 0.28 & 90 & 7.00 \\
P-n50-k7 & 50 & 1.6 & 1278.008 & 0.8 & 1277.174 & 118.6 & 0.07 & 83 & 0.18 \\
P-n50-k7 & 38 & 1.9 & 1536.758 & 1.0 & 1533.340 & 19.3 & 0.22 & 71 & $-$2.11 \\
P-n50-k8 & 60 & 1.3 & 1165.777 & 0.3 & 1165.761 & 2.0 & 0.00 & 77 & $-$0.02 \\
P-n50-k8 & 40 & 1.6 & 1540.974 & 0.2 & 1539.982 & 2.2 & 0.06 & 67 & $-$0.30 \\
P-n50-k8 & 30 & 1.9 & 1919.737 & 0.1 & 1917.504 & 1.8 & 0.12 & 65 & 1.01 \\
P-n50-k10 & 33 & 1.6 & 1798.096 & 0.0 & 1798.094 & 0.1 & 0.00 & 58 & 0.00 \\
P-n50-k10 & 25 & 1.9 & 2262.481 & 0.0 & 2262.006 & 0.1 & 0.02 & 57 & 0.61 \\
Avg &  &  &  & 0.8 &  & 29.0 & 0.10 & 70 & 1.80 \\
\bottomrule
\end{tabular}
\captionsetup{justification=centering}
\caption*{$^a$ Solution cost improvement. $^b$ Proportion of solution routes that are different from the routes observed in the optimal restocking solution. $^c$ Late service improvement.}
\endgroup
\end{table}

\begin{figure}[t]
\centering
	{\includegraphics[trim={45 45 45 45},clip,width=300pt]{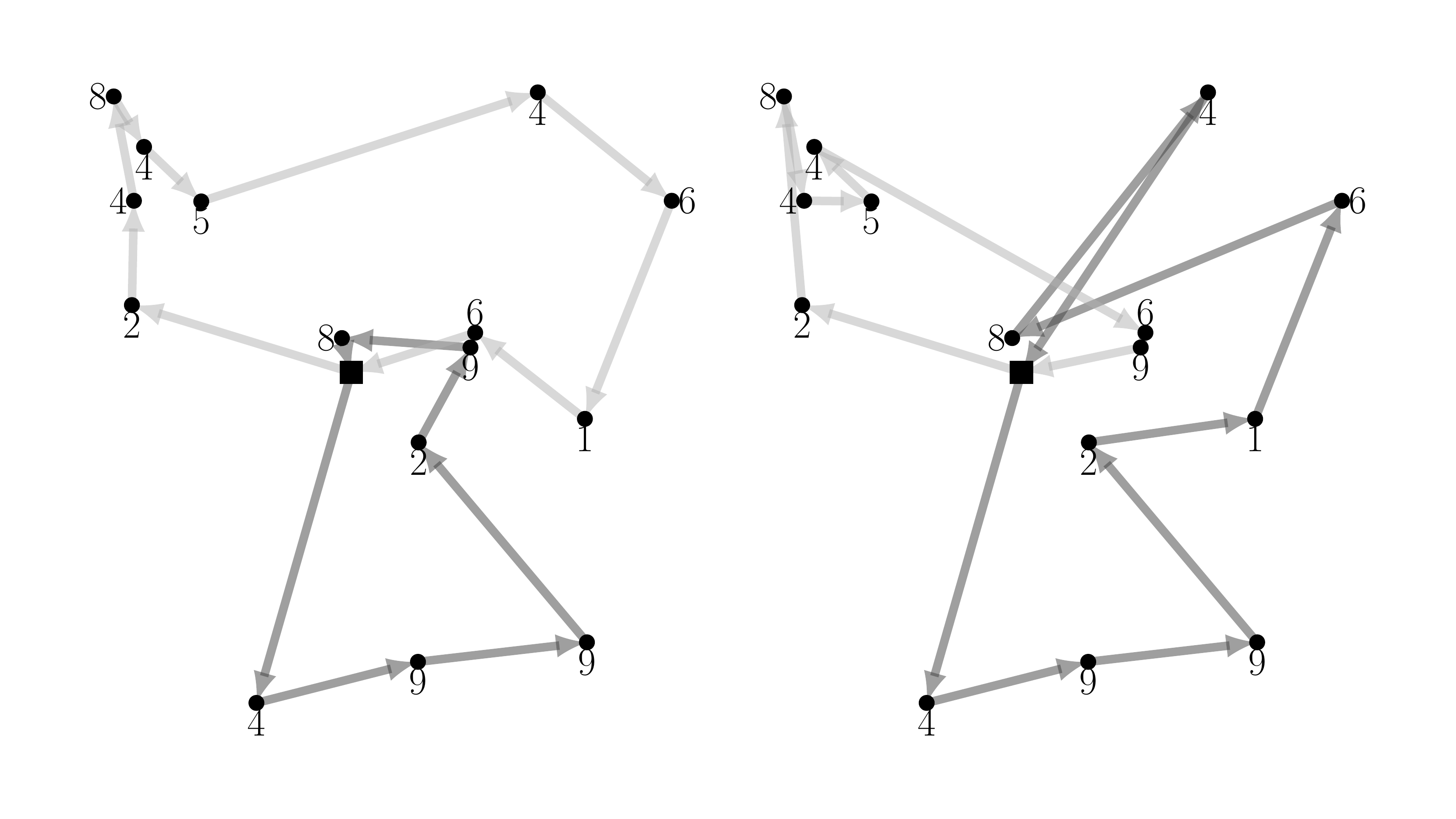}}
	{}
		\caption{Optimal planned routes under optimal restocking (left) and switch (right) policies. The optimal solution under the switch policy is $1.04\%$ less costly.\label{fig:solcomp}}
\end{figure}

More often than not, optimal solutions under the switch and optimal restocking policies display the same assignment of customers to vehicles, but differ in terms of the planned sequences along which customers are visited. Figure \ref{fig:solcomp} illustrates one example in which this structural similarity does not hold. In this instance with vehicle capacity $Q=23$, load factor $f=1.9$ and low variability demands, the optimal solution under the switch policy is $1.04\%$ less costly and is achieved by changing customer-vehicle assignments and planned sequences. Note also that the optimal solution under the switch policy exhibits planned sequences that do not correspond to optimal Hamiltonian cycles. This highlights the difficulty of the VRPSD under complex recourse policies when compared to the (deterministic) capacitated VRP: in the VRPSD, even when the customer-vehicle assignment is defined (as in the single-vehicle case), the problem of determining optimal planned sequences remains intractable except in very simple instances \citep{florio2020optimal}.

\begin{figure}[h]
\centering
	{\includegraphics[trim={55 155 55 125},clip,width=\textwidth]{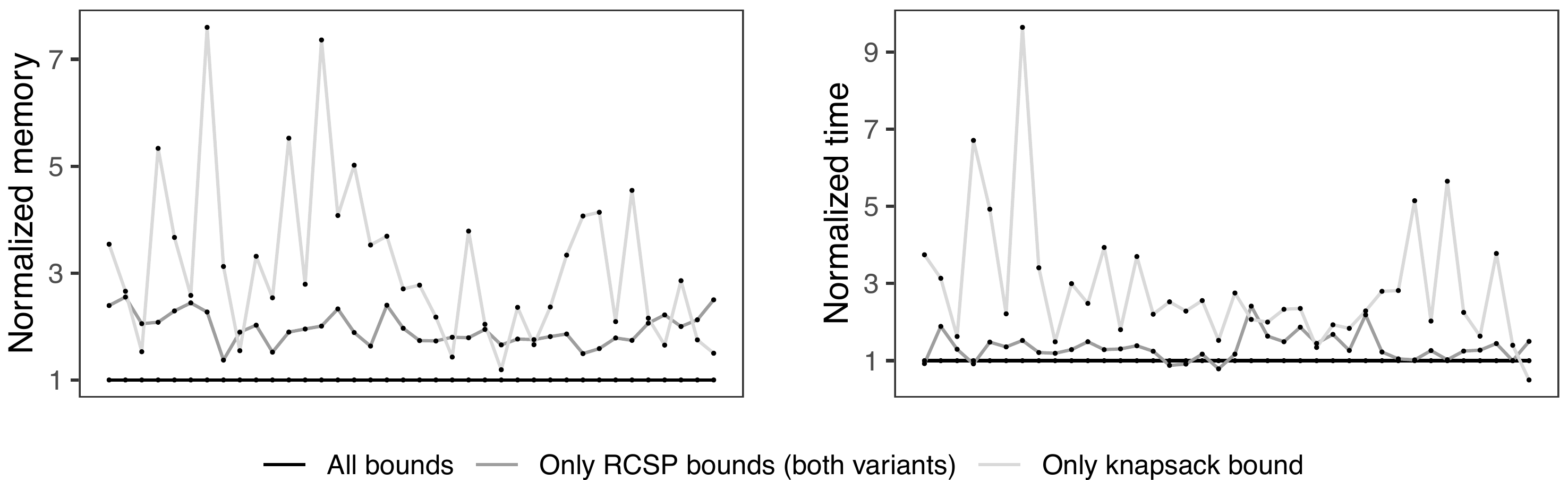}}
	{}
	\vspace{-9pt}

		\caption{Normalized memory and runtime requirements for solving the linear relaxation of 38 instances.\label{fig:bounds}}
\end{figure}

Finally, we investigate how effectively completion bounds can reduce the resource (memory and time) requirements of our algorithm. We focus on the resources needed to solve the RMP when only capacity cuts are separated. When all bounds are enabled, a lower bound can be computed in 58 out of the 63 instances with Poisson distributed demands. This number drops to 53 when only RCSP bounds (both variants) are enabled, and to 38 when only the knapsack bound is active. The plots in Figure \ref{fig:bounds} show the memory and runtime requirements of the algorithm over the 38 instances that can be solved with only the knapsack bound. Resource amounts are normalized, that is, they are expressed as a ratio of the amount required when all bounds are active. As shown, even though RCSP-based bounds are relatively more effective, both bounds, in combination, lead to considerably less memory usage. Runtime requirements decrease only slightly when the knapsack bound is also active (compared to only RCSP bounds), since the savings accrued by discarding more labels are partially offset by the additional runtime needed to solve a large number of knapsack problems.

\vspace{3pt}

\section{Conclusions} \label{sec:conclusions}
In this paper, we considered the VRPSD under a partial reoptimization recourse policy known as the switch policy. According to this recourse policy, the order of adjacent customers along an a priori sequence can be switched during route execution. By developing a branch-cut-and-price algorithm for the VRPSD under the switch policy, we were able to assess the value of reoptimization and to provide initial answers to a still open question on the VRPSD, namely, how much savings can be realized by allowing the reoptimization of a priori sequences. When planned sequences are not optimized for the recourse policy, cost savings by allowing customer swapping may be significant. On the other hand, the benefit of partial reoptimization relative to optimal restocking is mostly marginal when planned sequences are optimized for the recourse policy.

The design of a branch-cut-and-price algorithm able to handle such a complex recourse policy required a departure from traditional techniques usually implemented in branch-and-price methods for vehicle routing. This also led to new solution strategies which can inspire further developments for other complex VRPSDs. Most distinctly, in our algorithm we opted for elementary route pricing, after recognizing the difficulty of deriving effective dominance rules when labels must store complete value functions. In this setting, completion bounds of varying complexity are employed to control the combinatorial growth of labels when solving the subproblem. The proposed algorithm finds its limits in instances with many customers per route, as common in branch-and-price.

Besides measuring the value of reoptimization, extensive numerical experiments demonstrated other benefits of the switch policy such as reductions in the failure rates and in the number of customers served after replenishment. The results also show that a posteriori sequences change structurally depending on the restocking policy, and that the switch policy often employs customer reordering to reduce traveling costs. When considering optimized multi-vehicle plans, (partial) reoptimization brings modest cost savings in the instances solved, which are characterized by a relatively small number of customers per route. In this context, an interesting research perspective is to identify and solve instance classes where the switch recourse policy leads to larger savings. To this end, deriving tighter completion bounds could be instrumental in enabling the solution of a broader class of instances.

The literature on exact algorithms for the VRPSD is evolving towards solving the problem under increasingly complex recourse policies. Our results show lower cost savings (on average) when moving from optimal restocking to the switch policy, when compared to moving from detour-to-depot to optimal restocking. Therefore, future research should also focus on recourse strategies that go beyond intra-route reoptimization, such as, for example, the cooperative recourse strategies proposed by \cite{AkErera2007} and \cite{Zhu_2014}. While incorporating such complex policies within exact solution frameworks is certainly challenging, we expect that leveraging capacity pooling to prevent replenishment trips is likely to promote larger savings than yet more advanced intra-route reoptimization policies.

%\ACKNOWLEDGMENT{The authors thank the associate editor and two reviewers for many careful and insightful comments that helped to improve significantly this paper.}

\bibliography{SVRPRefs} % if more than one, comma separated

%\clearpage
%
%\setcounter{page}{1}
\appendix

\section{Detailed Branch-Cut-and-Price Results}

\setcounter{table}{0}
\renewcommand{\thetable}{A\arabic{table}}
Tables \ref{tab:bcpall13}-\ref{tab:bcpall19} present the complete results obtained by the branch-cut-and-price algorithm on 189 VRPSD instances. Each instance is experimented under load factors $f\in\{1.3,1.6,1.9\}$ and three probability distributions (column ``Dt''): Poisson (P), negative binomial with a variance-to-mean ratio of 2 (N), or binomial with a variance-to-mean ratio of 0.5 (B). The best solution found and the runtime (in minutes) are presented in columns ``Best'' and ``T'', respectively. Finally, columns ``Cols'', ``$\text{c}_{\textsf{av}}$'', ``$\text{c}_{\textsf{mx}}$'', ``CC'', ``RC'', ``LB'', ``T$_\text{LB}$'' and ``BB'' report the number of columns generated, the average number of customers per route in a solution, the maximum number of customers in a column generated, the number of strengthened capacity cuts separated, the number of subset row cuts separated, the lower bound obtained at the root node, the runtime (in minutes) for obtaining the root node lower bound, and the number of nodes explored throughout branch-and-bound, respectively.

\footnotesize

\begin{center}
\begin{longtable}{lrcrrrrrrrrrrr}
\caption{Branch-Cut-and-Price Results ($f=1.3$) \label{tab:bcpall13}} \\
\toprule
Instance & $Q$ & Dt & Best & Gap & T & Cols & $\text{c}_{\textsf{av}}$ & $\text{c}_{\textsf{mx}}$ & CC & RC & LB & T$_\text{LB}$ & BB \\
\midrule
\endfirsthead
\toprule
Instance & $Q$ & Dt & Best & Gap & T & Cols & $\text{c}_{\textsf{av}}$ & $\text{c}_{\textsf{mx}}$ & CC & RC & LB & T$_\text{LB}$ & BB \\
\midrule
\endhead
\midrule
\multicolumn{14}{r}{{Continued on next page}} \\
\bottomrule
\endfoot
\bottomrule
\endlastfoot
A-n32-k5 & 50 & P & 1538.008 & 0 & 80.4 & 2033 & 4.4 & 10 & 11 & 14 & 1538.008 & 80.4 & 1 \\
A-n32-k5 & 50 & N & 1551.732 & 0 & 58.7 & 1833 & 4.4 & 9 & 11 & 4 & 1551.732 & 58.7 & 1 \\
A-n32-k5 & 50 & B & 1523.919 & 1.169 & 120.0 & 2331 & 4.4 & 9 & 16 & 22 & 1505.870 & 55.3 & 5 \\
A-n33-k5 & 50 & P & 1223.964 & 0 & 30.4 & 2102 & 4.6 & 8 & 9 & 10 & 1223.584 & 5.2 & 15 \\
A-n33-k5 & 50 & N & 1231.986 & 0 & 2.4 & 1722 & 4.6 & 9 & 5 & 0 & 1231.986 & 2.4 & 1 \\
A-n33-k5 & 50 & B & 1199.517 & 0 & 20.5 & 2440 & 4.6 & 8 & 7 & 2 & 1199.496 & 1.8 & 19 \\
A-n34-k5 & 50 & P & 1460.946 & 0.014 & 120.0 & 2247 & 4.1 & 8 & 8 & 17 & 1460.310 & 27.0 & 21 \\
A-n34-k5 & 50 & N & 1470.423 & 0 & 73.7 & 1862 & 4.1 & 8 & 18 & 16 & 1470.164 & 44.6 & 3 \\
A-n34-k5 & 50 & B & 1432.669 & 0.209 & 120.0 & 2470 & 4.1 & 8 & 10 & 21 & 1428.946 & 11.4 & 33 \\
A-n36-k5 & 50 & P & 1497.431 & n/a & 120.0 & 2481 & 5.0 & 10 & 8 & 10 & n/a & n/a & 1 \\
A-n36-k5 & 50 & N & n/a & n/a & 120.0 & 2542 & 5.0 & 10 & 14 & 10 & n/a & n/a & 1 \\
A-n36-k5 & 50 & B & 1499.889 & n/a & 120.0 & 2439 & 5.0 & 9 & 5 & 20 & n/a & n/a & 1 \\
A-n37-k5 & 50 & P & 1102.222 & n/a & 120.0 & 3071 & 5.1 & 13 & 7 & 8 & n/a & n/a & 1 \\
A-n37-k5 & 50 & N & 1099.681 & 0 & 73.8 & 2803 & 5.1 & 14 & 6 & 0 & 1099.681 & 73.8 & 1 \\
A-n37-k5 & 50 & B & 1092.598 & n/a & 120.0 & 3220 & 5.1 & 14 & 10 & 16 & n/a & n/a & 1 \\
A-n38-k5 & 50 & P & 1331.285 & 0.008 & 120.0 & 2686 & 4.6 & 9 & 16 & 8 & 1331.180 & 43.4 & 5 \\
A-n38-k5 & 50 & N & 1335.740 & 0 & 31.3 & 2217 & 4.6 & 9 & 10 & 0 & 1335.740 & 31.3 & 1 \\
A-n38-k5 & 50 & B & 1321.567 & 0.369 & 120.0 & 2758 & 4.6 & 9 & 22 & 20 & 1316.288 & 29.4 & 7 \\
A-n39-k5 & 50 & P & n/a & n/a & 120.0 & 2671 & 4.8 & 11 & 10 & 16 & n/a & n/a & 1 \\
A-n39-k5 & 50 & N & 1566.337 & n/a & 120.0 & 2369 & 4.8 & 11 & 9 & 10 & n/a & n/a & 1 \\
A-n39-k5 & 50 & B & n/a & n/a & 120.0 & 2418 & 4.8 & 10 & 7 & 10 & n/a & n/a & 1 \\
A-n44-k6 & 50 & P & n/a & n/a & 120.0 & 3255 & 4.8 & 9 & 24 & 16 & n/a & n/a & 1 \\
A-n44-k6 & 50 & N & n/a & n/a & 120.0 & 2829 & 4.8 & 9 & 20 & 6 & n/a & n/a & 1 \\
A-n44-k6 & 50 & B & 1796.043 & 1.358 & 120.0 & 2896 & 4.8 & 9 & 15 & 18 & 1771.975 & 86.1 & 3 \\
A-n45-k6 & 50 & P & 1818.520 & 0.006 & 120.0 & 2651 & 4.4 & 10 & 23 & 4 & 1818.415 & 72.9 & 3 \\
A-n45-k6 & 50 & N & n/a & n/a & 120.0 & 2843 & 4.4 & 8 & 29 & 10 & n/a & n/a & 1 \\
A-n45-k6 & 50 & B & 1795.692 & 0 & 63.4 & 2580 & 4.4 & 9 & 20 & 10 & 1795.692 & 63.4 & 1 \\
A-n45-k7 & 50 & P & 2454.019 & n/a & 120.0 & 2490 & 4.4 & 8 & 11 & 8 & n/a & n/a & 1 \\
A-n45-k7 & 50 & N & 2476.503 & n/a & 120.0 & 2663 & 4.4 & 8 & 17 & 8 & n/a & n/a & 1 \\
A-n45-k7 & 50 & B & 2405.556 & n/a & 120.0 & 2638 & 4.4 & 9 & 11 & 10 & n/a & n/a & 1 \\
A-n46-k7 & 50 & P & 1853.487 & n/a & 120.0 & 2901 & 4.5 & 10 & 10 & 14 & n/a & n/a & 1 \\
A-n46-k7 & 50 & N & n/a & n/a & 120.0 & 2793 & 4.5 & 9 & 18 & 10 & n/a & n/a & 1 \\
A-n46-k7 & 50 & B & n/a & n/a & 120.0 & 2952 & 4.5 & 10 & 7 & 12 & n/a & n/a & 1 \\
A-n48-k7 & 50 & P & n/a & n/a & 120.0 & 3176 & 4.7 & 9 & 16 & 14 & n/a & n/a & 1 \\
A-n48-k7 & 50 & N & n/a & n/a & 120.0 & 3183 & 4.7 & 9 & 18 & 6 & n/a & n/a & 1 \\
A-n48-k7 & 50 & B & n/a & n/a & 120.0 & 3082 & 4.7 & 9 & 8 & 18 & n/a & n/a & 1 \\
P-n19-k2 & 80 & P & 370.006 & 0 & 3.6 & 1128 & 6.0 & 8 & 0 & 2 & 370.006 & 3.6 & 1 \\
P-n19-k2 & 80 & N & 370.659 & 0 & 6.8 & 1142 & 6.0 & 9 & 3 & 2 & 370.659 & 6.8 & 1 \\
P-n19-k2 & 80 & B & 367.393 & 0 & 1.7 & 1018 & 6.0 & 9 & 1 & 0 & 367.393 & 1.7 & 1 \\
P-n20-k2 & 80 & P & 375.066 & 0 & 9.4 & 1319 & 6.3 & 9 & 0 & 2 & 375.066 & 9.4 & 1 \\
P-n20-k2 & 80 & N & 374.315 & 0 & 17.6 & 1294 & 6.3 & 9 & 1 & 2 & 374.315 & 17.6 & 1 \\
P-n20-k2 & 80 & B & 373.747 & 0 & 11.0 & 1336 & 6.3 & 10 & 2 & 4 & 373.747 & 11.0 & 1 \\
P-n21-k2 & 80 & P & 363.013 & 0 & 42.6 & 1683 & 6.7 & 10 & 5 & 8 & 363.013 & 42.6 & 1 \\
P-n21-k2 & 80 & N & 363.029 & 0 & 49.2 & 1557 & 6.7 & 10 & 5 & 2 & 363.029 & 49.2 & 1 \\
P-n21-k2 & 80 & B & 361.718 & 0.147 & 120.0 & 2016 & 6.7 & 10 & 3 & 10 & 361.165 & 31.5 & 7 \\
P-n22-k2 & 80 & P & 375.197 & 0 & 43.1 & 1698 & 7.0 & 10 & 0 & 0 & 375.197 & 43.1 & 1 \\
P-n22-k2 & 80 & N & 375.548 & 0 & 62.0 & 1671 & 7.0 & 10 & 0 & 0 & 375.548 & 62.0 & 1 \\
P-n22-k2 & 80 & B & 373.750 & 0 & 39.5 & 1770 & 7.0 & 10 & 1 & 0 & 373.750 & 39.5 & 1 \\
P-n40-k5 & 70 & P & 769.656 & n/a & 120.0 & 2975 & 5.6 & 10 & 8 & 10 & n/a & n/a & 1 \\
P-n40-k5 & 70 & N & 771.150 & n/a & 120.0 & 2972 & 5.6 & 9 & 14 & 6 & n/a & n/a & 1 \\
P-n40-k5 & 70 & B & 768.467 & 0.756 & 120.0 & 3159 & 5.6 & 10 & 9 & 14 & 762.698 & 93.4 & 3 \\
P-n45-k5 & 75 & P & n/a & n/a & 120.0 & 3534 & 5.5 & 11 & 7 & 10 & n/a & n/a & 1 \\
P-n45-k5 & 75 & N & n/a & n/a & 120.0 & 3444 & 5.5 & 10 & 10 & 6 & n/a & n/a & 1 \\
P-n45-k5 & 75 & B & n/a & n/a & 120.0 & 3624 & 5.5 & 10 & 3 & 10 & n/a & n/a & 1 \\
P-n50-k7 & 75 & P & 1002.557 & 0.764 & 120.0 & 3277 & 4.9 & 8 & 12 & 18 & 994.955 & 62.6 & 3 \\
P-n50-k7 & 75 & N & 1002.747 & 0.382 & 120.0 & 3254 & 4.9 & 8 & 15 & 14 & 998.935 & 116.4 & 3 \\
P-n50-k7 & 75 & B & 991.272 & 0.311 & 120.0 & 3243 & 4.9 & 8 & 15 & 14 & 988.122 & 58.3 & 7 \\
P-n50-k8 & 60 & P & 1165.761 & 0 & 2.0 & 1923 & 3.8 & 7 & 6 & 0 & 1165.761 & 2.0 & 1 \\
P-n50-k8 & 60 & N & 1174.268 & 0 & 8.7 & 1905 & 3.8 & 7 & 23 & 6 & 1174.268 & 8.7 & 1 \\
P-n50-k8 & 60 & B & 1150.529 & 0 & 2.0 & 1935 & 3.8 & 7 & 7 & 2 & 1150.529 & 2.0 & 1 \\
P-n50-k10 & 50 & P & 1380.337 & 0.123 & 120.0 & 3844 & 3.3 & 6 & 90 & 57 & 1375.903 & 1.0 & 351 \\
P-n50-k10 & 50 & N & 1380.698 & 0 & 25.3 & 2277 & 3.3 & 6 & 45 & 28 & 1378.725 & 1.7 & 45 \\
P-n50-k10 & 50 & B & 1368.938 & 0.395 & 120.0 & 4105 & 3.3 & 6 & 170 & 75 & 1361.515 & 0.7 & 433 \\
\end{longtable}
\end{center}

\begin{center}
\begin{longtable}{lrcrrrrrrrrrrr}
\caption{Branch-Cut-and-Price Results ($f=1.6$) \label{tab:bcpall16}} \\
\toprule
Instance & $Q$ & Dt & Best & Gap & T & Cols & $\text{c}_{\textsf{av}}$ & $\text{c}_{\textsf{mx}}$ & CC & RC & LB & T$_\text{LB}$ & BB \\
\midrule
\endfirsthead
\toprule
Instance & $Q$ & Dt & Best & Gap & T & Cols & $\text{c}_{\textsf{av}}$ & $\text{c}_{\textsf{mx}}$ & CC & RC & LB & T$_\text{LB}$ & BB \\
\midrule
\endhead
\midrule
\multicolumn{14}{r}{{Continued on next page}} \\
\bottomrule
\endfoot
\bottomrule
\endlastfoot
A-n32-k5 & 33 & P & 2110.029 & 0.249 & 120.0 & 1842 & 3.9 & 8 & 33 & 20 & 2104.668 & 26.7 & 11 \\
A-n32-k5 & 33 & N & 2146.480 & 0.154 & 120.0 & 1884 & 3.9 & 7 & 22 & 12 & 2143.180 & 29.3 & 9 \\
A-n32-k5 & 33 & B & 2076.346 & 0.301 & 120.0 & 1744 & 3.9 & 9 & 20 & 18 & 2068.988 & 31.3 & 11 \\
A-n33-k5 & 33 & P & 1633.791 & 0.004 & 120.0 & 3479 & 3.6 & 7 & 56 & 38 & 1627.303 & 1.5 & 311 \\
A-n33-k5 & 33 & N & 1665.618 & 0 & 63.2 & 2062 & 3.6 & 7 & 38 & 29 & 1661.507 & 3.6 & 63 \\
A-n33-k5 & 33 & B & 1605.823 & 0 & 29.4 & 1921 & 3.6 & 7 & 27 & 18 & 1602.243 & 1.2 & 63 \\
A-n34-k5 & 33 & P & 1961.076 & 0 & 59.3 & 1613 & 3.7 & 8 & 26 & 16 & 1959.205 & 5.9 & 19 \\
A-n34-k5 & 33 & N & 1996.306 & 0 & 7.6 & 1218 & 3.7 & 7 & 18 & 4 & 1996.306 & 7.6 & 1 \\
A-n34-k5 & 33 & B & 1926.803 & 0 & 56.5 & 1911 & 3.7 & 7 & 28 & 12 & 1925.735 & 6.5 & 25 \\
A-n36-k5 & 33 & P & 2026.418 & 0.546 & 120.0 & 1855 & 3.9 & 9 & 27 & 22 & 2015.167 & 41.3 & 5 \\
A-n36-k5 & 33 & N & 2062.656 & 0.573 & 120.0 & 1786 & 3.9 & 9 & 20 & 12 & 2050.896 & 56.9 & 3 \\
A-n36-k5 & 33 & B & 1998.868 & 0.718 & 120.0 & 1910 & 3.9 & 9 & 19 & 16 & 1984.419 & 38.6 & 5 \\
A-n37-k5 & 33 & P & n/a & n/a & 120.0 & 2304 & 4.5 & 13 & 19 & 12 & n/a & n/a & 1 \\
A-n37-k5 & 33 & N & n/a & n/a & 120.0 & 2393 & 4.5 & 13 & 23 & 14 & n/a & n/a & 1 \\
A-n37-k5 & 33 & B & 1418.943 & 0.343 & 120.0 & 2425 & 4.5 & 11 & 14 & 12 & 1414.089 & 92.6 & 3 \\
A-n38-k5 & 33 & P & 1772.391 & 0.063 & 120.0 & 1958 & 3.7 & 7 & 29 & 12 & 1770.318 & 20.9 & 15 \\
A-n38-k5 & 33 & N & 1796.386 & 0 & 22.8 & 1498 & 3.7 & 7 & 22 & 6 & 1796.386 & 22.8 & 1 \\
A-n38-k5 & 33 & B & 1736.760 & 0 & 89.3 & 2008 & 3.7 & 8 & 24 & 10 & 1734.702 & 23.3 & 9 \\
A-n39-k5 & 33 & P & 2083.670 & 0.734 & 120.0 & 2178 & 4.2 & 9 & 21 & 12 & 2068.489 & 75.0 & 3 \\
A-n39-k5 & 33 & N & 2131.285 & 0.361 & 120.0 & 2059 & 4.2 & 9 & 19 & 20 & 2123.628 & 103.2 & 3 \\
A-n39-k5 & 33 & B & 2031.663 & 0.342 & 120.0 & 2116 & 4.2 & 9 & 15 & 10 & 2024.609 & 49.4 & 5 \\
A-n44-k6 & 33 & P & 2399.516 & 0.348 & 120.0 & 2232 & 3.9 & 8 & 24 & 16 & 2390.054 & 52.3 & 7 \\
A-n44-k6 & 33 & N & 2448.279 & 0.452 & 120.0 & 2148 & 3.9 & 8 & 32 & 18 & 2437.265 & 74.8 & 3 \\
A-n44-k6 & 33 & B & 2362.595 & 0.399 & 120.0 & 2218 & 3.9 & 8 & 29 & 22 & 2352.341 & 43.4 & 5 \\
A-n45-k6 & 33 & P & 2439.659 & 0.145 & 120.0 & 2180 & 3.7 & 7 & 37 & 16 & 2435.894 & 38.3 & 7 \\
A-n45-k6 & 33 & N & 2482.809 & 0.082 & 120.0 & 2005 & 3.7 & 8 & 46 & 14 & 2480.779 & 57.7 & 3 \\
A-n45-k6 & 33 & B & 2404.644 & 0.247 & 120.0 & 2167 & 3.7 & 8 & 32 & 22 & 2396.256 & 25.1 & 7 \\
A-n45-k7 & 33 & P & 3297.135 & 0.224 & 120.0 & 1881 & 3.4 & 7 & 40 & 8 & 3289.760 & 45.1 & 5 \\
A-n45-k7 & 33 & N & 3371.065 & 0.209 & 120.0 & 1712 & 3.4 & 7 & 37 & 10 & 3364.037 & 63.0 & 3 \\
A-n45-k7 & 33 & B & 3228.300 & 0.117 & 120.0 & 1945 & 3.4 & 7 & 35 & 14 & 3224.370 & 43.1 & 5 \\
A-n46-k7 & 33 & P & 2488.238 & 0.951 & 120.0 & 2210 & 3.8 & 7 & 29 & 16 & 2464.792 & 75.4 & 3 \\
A-n46-k7 & 33 & N & 2526.362 & n/a & 120.0 & 1983 & 3.8 & 8 & 32 & 8 & n/a & n/a & 1 \\
A-n46-k7 & 33 & B & 2441.972 & 0.687 & 120.0 & 2109 & 3.8 & 7 & 14 & 12 & 2425.318 & 70.6 & 3 \\
A-n48-k7 & 33 & P & 3011.586 & n/a & 120.0 & 2296 & 3.9 & 8 & 29 & 10 & n/a & n/a & 1 \\
A-n48-k7 & 33 & N & 3064.337 & n/a & 120.0 & 2240 & 3.9 & 7 & 29 & 10 & n/a & n/a & 1 \\
A-n48-k7 & 33 & B & 2958.475 & 0.272 & 120.0 & 2429 & 3.9 & 7 & 25 & 10 & 2950.443 & 102.7 & 3 \\
P-n19-k2 & 53 & P & 469.539 & 0 & 33.4 & 2146 & 4.5 & 7 & 10 & 33 & 467.511 & 2.0 & 61 \\
P-n19-k2 & 53 & N & 473.578 & 0 & 10.5 & 1067 & 4.5 & 7 & 6 & 14 & 473.452 & 2.5 & 13 \\
P-n19-k2 & 53 & B & 467.404 & 0.216 & 120.0 & 6546 & 4.5 & 7 & 30 & 126 & 462.936 & 2.5 & 431 \\
P-n20-k2 & 53 & P & 480.387 & 0.476 & 120.0 & 2267 & 4.8 & 8 & 23 & 68 & 476.361 & 6.6 & 57 \\
P-n20-k2 & 53 & N & 483.374 & 0 & 101.7 & 1738 & 4.8 & 8 & 11 & 20 & 482.426 & 10.9 & 37 \\
P-n20-k2 & 53 & B & 477.091 & 0.936 & 120.0 & 2753 & 4.8 & 8 & 16 & 80 & 469.918 & 2.7 & 67 \\
P-n21-k2 & 53 & P & 466.729 & 0 & 8.8 & 1083 & 5.0 & 9 & 5 & 2 & 466.729 & 8.8 & 1 \\
P-n21-k2 & 53 & N & 474.804 & 0 & 12.7 & 1078 & 5.0 & 9 & 7 & 4 & 474.804 & 12.7 & 1 \\
P-n21-k2 & 53 & B & 462.180 & 0 & 10.4 & 1002 & 5.0 & 9 & 0 & 6 & 462.180 & 10.4 & 1 \\
P-n22-k2 & 53 & P & 487.682 & 0.896 & 120.0 & 1595 & 5.3 & 9 & 9 & 34 & 482.969 & 37.9 & 7 \\
P-n22-k2 & 53 & N & 490.748 & 0.098 & 120.0 & 1627 & 5.3 & 8 & 15 & 18 & 490.146 & 47.5 & 7 \\
P-n22-k2 & 53 & B & 481.293 & 0.81 & 120.0 & 1810 & 5.3 & 9 & 9 & 38 & 476.661 & 25.6 & 9 \\
P-n40-k5 & 47 & P & 963.165 & 0.681 & 120.0 & 2443 & 4.3 & 9 & 27 & 12 & 956.626 & 37.7 & 5 \\
P-n40-k5 & 47 & N & 976.540 & 0.011 & 120.0 & 2110 & 4.3 & 8 & 23 & 10 & 976.434 & 74.3 & 5 \\
P-n40-k5 & 47 & B & 946.874 & 0.371 & 120.0 & 2420 & 4.3 & 8 & 16 & 16 & 943.195 & 31.3 & 9 \\
P-n45-k5 & 50 & P & n/a & n/a & 120.0 & 2838 & 4.9 & 9 & 24 & 12 & n/a & n/a & 1 \\
P-n45-k5 & 50 & N & n/a & n/a & 120.0 & 2863 & 4.9 & 9 & 22 & 12 & n/a & n/a & 1 \\
P-n45-k5 & 50 & B & 1055.302 & 0.792 & 120.0 & 2957 & 4.9 & 9 & 16 & 12 & 1047.011 & 106.0 & 3 \\
P-n50-k7 & 50 & P & 1277.174 & 0 & 118.6 & 2684 & 4.1 & 7 & 18 & 14 & 1275.704 & 9.3 & 35 \\
P-n50-k7 & 50 & N & 1303.486 & 0.379 & 120.0 & 2568 & 4.1 & 7 & 32 & 30 & 1297.481 & 10.3 & 15 \\
P-n50-k7 & 50 & B & 1263.822 & 0 & 15.3 & 2185 & 4.1 & 7 & 19 & 12 & 1263.822 & 15.3 & 1 \\
P-n50-k8 & 40 & P & 1539.982 & 0 & 2.2 & 1467 & 3.3 & 6 & 26 & 6 & 1539.147 & 0.7 & 5 \\
P-n50-k8 & 40 & N & 1575.684 & 0 & 11.2 & 1734 & 3.3 & 6 & 29 & 4 & 1575.649 & 0.9 & 25 \\
P-n50-k8 & 40 & B & 1518.116 & 0 & 1.7 & 1517 & 3.3 & 6 & 25 & 10 & 1518.116 & 1.7 & 1 \\
P-n50-k10 & 33 & P & 1798.094 & 0 & 0.1 & 897 & 2.6 & 5 & 25 & 2 & 1798.094 & 0.1 & 1 \\
P-n50-k10 & 33 & N & 1844.983 & 0 & 0.2 & 929 & 2.6 & 5 & 21 & 4 & 1844.983 & 0.2 & 1 \\
P-n50-k10 & 33 & B & 1756.539 & 0 & 5.6 & 1356 & 2.6 & 5 & 55 & 10 & 1754.896 & 0.1 & 83 \\
\end{longtable}
\end{center}

\begin{center}
\begin{longtable}{lrcrrrrrrrrrrr}
\caption{Branch-Cut-and-Price Results ($f=1.9$) \label{tab:bcpall19}} \\
\toprule
Instance & $Q$ & Dt & Best & Gap & T & Cols & $\text{c}_{\textsf{av}}$ & $\text{c}_{\textsf{mx}}$ & CC & RC & LB & T$_\text{LB}$ & BB \\
\midrule
\endfirsthead
\toprule
Instance & $Q$ & Dt & Best & Gap & T & Cols & $\text{c}_{\textsf{av}}$ & $\text{c}_{\textsf{mx}}$ & CC & RC & LB & T$_\text{LB}$ & BB \\
\midrule
\endhead
\midrule
\multicolumn{14}{r}{{Continued on next page}} \\
\bottomrule
\endfoot
\bottomrule
\endlastfoot
A-n32-k5 & 25 & P & 2615.913 & 0 & 39.5 & 1397 & 3.4 & 9 & 18 & 9 & 2615.115 & 13.1 & 5 \\
A-n32-k5 & 25 & N & 2661.403 & 0 & 12.1 & 1204 & 3.4 & 8 & 20 & 6 & 2661.403 & 12.1 & 1 \\
A-n32-k5 & 25 & B & 2552.970 & 0 & 91.7 & 1891 & 3.4 & 7 & 25 & 14 & 2543.686 & 8.8 & 43 \\
A-n33-k5 & 25 & P & 2018.620 & 0 & 22.4 & 1562 & 3.2 & 7 & 39 & 19 & 2014.464 & 1.3 & 61 \\
A-n33-k5 & 25 & N & 2058.124 & 0 & 24.9 & 1742 & 3.2 & 6 & 31 & 23 & 2052.389 & 0.8 & 67 \\
A-n33-k5 & 25 & B & 1958.271 & 0 & 0.7 & 1000 & 3.2 & 7 & 23 & 10 & 1958.271 & 0.7 & 1 \\
A-n34-k5 & 25 & P & 2429.258 & 0 & 46.7 & 1541 & 3.3 & 7 & 29 & 8 & 2427.511 & 1.7 & 33 \\
A-n34-k5 & 25 & N & 2469.894 & 0 & 17.2 & 1245 & 3.3 & 6 & 22 & 0 & 2469.802 & 1.1 & 13 \\
A-n34-k5 & 25 & B & 2373.021 & 0 & 4.4 & 1181 & 3.3 & 7 & 22 & 6 & 2371.759 & 0.9 & 5 \\
A-n36-k5 & 25 & P & 2472.101 & 0.018 & 120.0 & 1703 & 3.5 & 8 & 22 & 8 & 2471.661 & 23.4 & 11 \\
A-n36-k5 & 25 & N & 2510.435 & 0 & 23.9 & 1436 & 3.5 & 9 & 26 & 4 & 2510.435 & 23.9 & 1 \\
A-n36-k5 & 25 & B & 2408.964 & 0 & 14.3 & 1423 & 3.5 & 9 & 23 & 8 & 2408.964 & 14.3 & 1 \\
A-n37-k5 & 25 & P & n/a & n/a & 120.0 & 2120 & 4.0 & 12 & 25 & 4 & n/a & n/a & 1 \\
A-n37-k5 & 25 & N & n/a & n/a & 120.0 & 2058 & 4.0 & 12 & 29 & 2 & n/a & n/a & 1 \\
A-n37-k5 & 25 & B & n/a & n/a & 120.0 & 2228 & 4.0 & 12 & 16 & 2 & n/a & n/a & 1 \\
A-n38-k5 & 25 & P & 2171.219 & 0.092 & 120.0 & 2034 & 3.4 & 7 & 51 & 22 & 2167.404 & 13.3 & 25 \\
A-n38-k5 & 25 & N & 2211.356 & 0 & 9.9 & 1371 & 3.4 & 6 & 30 & 4 & 2211.356 & 9.9 & 1 \\
A-n38-k5 & 25 & B & 2111.597 & 0.171 & 120.0 & 2203 & 3.4 & 7 & 54 & 38 & 2104.858 & 6.5 & 41 \\
A-n39-k5 & 25 & P & 2548.391 & 0 & 42.6 & 1472 & 3.5 & 9 & 26 & 2 & 2548.391 & 42.6 & 1 \\
A-n39-k5 & 25 & N & 2600.107 & 0.137 & 120.0 & 1734 & 3.5 & 10 & 37 & 4 & 2596.548 & 47.5 & 3 \\
A-n39-k5 & 25 & B & 2477.731 & 0 & 92.5 & 1700 & 3.5 & 11 & 27 & 2 & 2476.688 & 19.1 & 5 \\
A-n44-k6 & 25 & P & 2917.078 & 0 & 49.1 & 1627 & 3.3 & 6 & 30 & 6 & 2916.907 & 15.1 & 5 \\
A-n44-k6 & 25 & N & 2982.691 & 0.079 & 120.0 & 1804 & 3.3 & 7 & 33 & 16 & 2979.273 & 29.5 & 13 \\
A-n44-k6 & 25 & B & 2825.973 & 0 & 55.6 & 1696 & 3.3 & 7 & 23 & 0 & 2825.814 & 4.6 & 15 \\
A-n45-k6 & 25 & P & 2997.605 & 0 & 8.9 & 1632 & 3.4 & 7 & 29 & 0 & 2997.605 & 8.9 & 1 \\
A-n45-k6 & 25 & N & 3069.523 & 0 & 19.9 & 1715 & 3.4 & 7 & 29 & 6 & 3069.523 & 19.9 & 1 \\
A-n45-k6 & 25 & B & 2914.347 & 0 & 6.2 & 1719 & 3.4 & 7 & 20 & 0 & 2914.347 & 6.2 & 1 \\
A-n45-k7 & 25 & P & 4099.811 & 0.011 & 120.0 & 1794 & 3.1 & 7 & 40 & 6 & 4098.182 & 19.9 & 11 \\
A-n45-k7 & 25 & N & 4177.747 & 0.002 & 120.0 & 1668 & 3.1 & 7 & 32 & 6 & 4177.514 & 25.5 & 11 \\
A-n45-k7 & 25 & B & 3989.004 & 0.053 & 120.0 & 1878 & 3.1 & 7 & 57 & 20 & 3981.357 & 14.2 & 15 \\
A-n46-k7 & 25 & P & 3073.507 & 0.028 & 120.0 & 1928 & 3.5 & 8 & 39 & 8 & 3072.639 & 58.6 & 5 \\
A-n46-k7 & 25 & N & 3133.525 & 0.263 & 120.0 & 1950 & 3.5 & 7 & 39 & 14 & 3125.305 & 67.3 & 3 \\
A-n46-k7 & 25 & B & 2987.472 & 0 & 33.9 & 1979 & 3.5 & 7 & 33 & 10 & 2987.472 & 33.9 & 1 \\
A-n48-k7 & 25 & P & 3691.833 & 0.113 & 120.0 & 1845 & 3.4 & 7 & 29 & 0 & 3691.383 & 44.0 & 3 \\
A-n48-k7 & 25 & N & 3773.098 & 0.101 & 120.0 & 1853 & 3.4 & 8 & 34 & 2 & 3769.290 & 55.6 & 3 \\
A-n48-k7 & 25 & B & 3585.691 & 0 & 33.2 & 1799 & 3.4 & 8 & 30 & 6 & 3585.691 & 33.2 & 1 \\
P-n19-k2 & 40 & P & 557.394 & 0 & 1.6 & 789 & 3.6 & 6 & 10 & 10 & 556.586 & 0.5 & 7 \\
P-n19-k2 & 40 & N & 569.175 & 0 & 0.8 & 584 & 3.6 & 7 & 6 & 6 & 569.175 & 0.8 & 1 \\
P-n19-k2 & 40 & B & 546.253 & 0 & 5.6 & 1066 & 3.6 & 6 & 8 & 20 & 544.483 & 0.7 & 31 \\
P-n20-k2 & 40 & P & 566.236 & 0 & 4.5 & 778 & 3.8 & 7 & 6 & 8 & 566.197 & 1.7 & 5 \\
P-n20-k2 & 40 & N & 579.043 & 0 & 2.2 & 716 & 3.8 & 7 & 4 & 4 & 579.043 & 2.2 & 1 \\
P-n20-k2 & 40 & B & 554.497 & 0 & 9.6 & 1012 & 3.8 & 7 & 12 & 14 & 553.419 & 2.3 & 13 \\
P-n21-k2 & 40 & P & 558.616 & 0 & 2.9 & 1009 & 5.0 & 8 & 3 & 0 & 558.616 & 2.9 & 1 \\
P-n21-k2 & 40 & N & 568.670 & 0 & 2.9 & 978 & 5.0 & 8 & 3 & 0 & 568.670 & 2.9 & 1 \\
P-n21-k2 & 40 & B & 547.592 & 0 & 1.9 & 1029 & 5.0 & 8 & 4 & 0 & 547.592 & 1.9 & 1 \\
P-n22-k2 & 40 & P & 573.906 & 0 & 5.0 & 874 & 4.2 & 8 & 6 & 2 & 573.906 & 5.0 & 1 \\
P-n22-k2 & 40 & N & 589.761 & 0 & 20.9 & 855 & 4.2 & 8 & 6 & 10 & 589.761 & 20.9 & 1 \\
P-n22-k2 & 40 & B & 561.078 & 0 & 11.8 & 898 & 4.2 & 8 & 4 & 11 & 561.078 & 11.8 & 1 \\
P-n40-k5 & 35 & P & 1172.855 & 0.004 & 120.0 & 1951 & 3.9 & 7 & 22 & 4 & 1172.814 & 17.5 & 17 \\
P-n40-k5 & 35 & N & 1202.193 & 0.218 & 120.0 & 1944 & 3.9 & 7 & 32 & 20 & 1200.798 & 43.8 & 5 \\
P-n40-k5 & 35 & B & 1142.283 & 0 & 29.0 & 1885 & 3.9 & 8 & 23 & 2 & 1142.197 & 13.9 & 3 \\
P-n45-k5 & 38 & P & 1266.367 & 0 & 103.9 & 2548 & 4.4 & 8 & 30 & 6 & 1266.258 & 49.9 & 3 \\
P-n45-k5 & 38 & N & 1290.375 & 0 & 76.7 & 2394 & 4.4 & 8 & 27 & 4 & 1290.375 & 76.7 & 1 \\
P-n45-k5 & 38 & B & 1236.093 & 0 & 50.4 & 2356 & 4.4 & 8 & 22 & 8 & 1236.093 & 50.4 & 1 \\
P-n50-k7 & 38 & P & 1533.340 & 0 & 19.3 & 1881 & 3.5 & 6 & 28 & 6 & 1532.879 & 5.3 & 7 \\
P-n50-k7 & 38 & N & 1571.214 & 0.052 & 120.0 & 2269 & 3.5 & 6 & 48 & 14 & 1568.616 & 8.0 & 23 \\
P-n50-k7 & 38 & B & 1495.087 & 0 & 38.8 & 2010 & 3.5 & 6 & 20 & 14 & 1494.493 & 7.3 & 15 \\
P-n50-k8 & 30 & P & 1917.504 & 0 & 1.8 & 1198 & 2.9 & 5 & 33 & 2 & 1917.254 & 0.3 & 11 \\
P-n50-k8 & 30 & N & 1946.668 & 0 & 0.5 & 1100 & 2.9 & 5 & 28 & 4 & 1946.668 & 0.5 & 1 \\
P-n50-k8 & 30 & B & 1874.369 & 0 & 6.1 & 1430 & 2.9 & 5 & 41 & 4 & 1873.381 & 0.3 & 47 \\
P-n50-k10 & 25 & P & 2262.006 & 0 & 0.1 & 710 & 2.3 & 5 & 27 & 4 & 2261.745 & 0.1 & 3 \\
P-n50-k10 & 25 & N & 2300.052 & 0 & 0.1 & 673 & 2.3 & 5 & 25 & 6 & 2300.052 & 0.1 & 1 \\
P-n50-k10 & 25 & B & 2209.437 & 0 & 0.1 & 676 & 2.3 & 5 & 20 & 4 & 2209.437 & 0.1 & 1 \\
\end{longtable}
\end{center}

\end{document}